\documentclass{article}
\usepackage{amssymb, amsmath,epsfig}
\usepackage{color}

\newcommand{\mc}[1]{{}}

\newcommand \mf {{\cal MF}}
\newcommand \ml {{\cal ML}}
\newcommand \pmf {{\cal PMF}}

\newcommand \La {\Lambda}

\newcommand{\qonetg}{{\cal Q}^1 {\cal T}_g}
\newcommand{\qonemg}{{\cal Q}^1 {\cal M}_g}
\newcommand{\qstartg}{{\cal Q}^* {\cal T}_g}

\newcommand{\omegaonetg}{\Omega^1 {\cal T}_g}

\newcommand {\mm} {{\bf m}}

\newcommand \F {{\cal F}}

\DeclareMathRadical{\sqrtsign}{symbols}{"70}{largesymbols}{"70}
\newcommand{\bb}{\mathbb}

\newcommand{\cx}{{\bb C}}

\newcommand{\integers}{{\bb Z}}

\newcommand{\reals}{{\bb R}}

\newlength{\figboxwidth}             
\newcommand{\makefig}[3]{
        \begin{figure}[htb]
        \refstepcounter{figure}
        \label{#2}
        \begin{center}
                #3~\\
                \smallskip
                Figure \thefigure.  #1
        \end{center}
        \medskip
        \end{figure}
}
\newcommand{\makefignocenter}[3]{
        \begin{figure}[htb]
        \refstepcounter{figure}
        \label{#2}
        \begin{center}
                #3~\\
        \end{center}
                \smallskip

                Figure \thefigure.  #1
        \medskip
        \end{figure}
}
\newcommand{\appendixmode}{
        \setcounter{section}{0}
        \renewcommand{\thesection}{\Alph{section}}
}

\renewcommand{\bold}[1]{\medskip \noindent {\bf #1 }\nopagebreak}

\newcommand{\qed}{\hfill $\Box$}

\newcommand{\tensor}{\otimes}

\newcommand{\cross}{\times}
\renewcommand{\Re}{\mbox{Re}\:}
\renewcommand{\Im}{\mbox{Im}\:}

\newcommand{\st}{\;\: : \;\:}         

\newcommand{\zed}{\integers}

\def\@ifundefined#1#2#3%
  {\expandafter\ifx\csname#1\endcsname\relax#2\else#3\fi}

\@ifundefined{theoremstyle}{
}{
\theoremstyle{plain} 
}
\newtheorem{theorem}{Theorem}[section]

\newtheorem{proposition}[theorem]{Proposition}
\newtheorem{lemma}[theorem]{Lemma}

\@ifundefined{theoremstyle}{
}{
\theoremstyle{definition} 
}

\newcommand{\cA}{{\cal A}}
\newcommand{\cB}{{\cal B}}

\newcommand{\cF}{{\cal F}}
\newcommand{\cG}{{\cal G}}
\newcommand{\cH}{{\cal H}}

\newcommand{\cK}{{\cal K}}

\newcommand{\cM}{{\cal M}}

\newcommand{\cP}{{\cal P}}
\newcommand{\cQ}{{\cal Q}}

\newcommand{\cS}{{\cal S}}
\newcommand{\cT}{{\cal T}}
\newcommand{\cU}{{\cal U}}
\newcommand{\cV}{{\cal V}}

\mathchardef\GG="321D
\catcode`~=11 
\newcommand{\urltilde}{\kern -.15em\lower .7ex\hbox{~}\kern .04em}
\catcode`~=13 

\newtheorem{theo}{Theorem}

\newtheorem{coro}[theo]{Corollary}

\newcommand{\te}{\mathcal{T}}
\newcommand{\vo}{\operatorname{Vol}}

\newcommand{\MF}{\mathcal{MF}}

\newcommand{\tw}{\operatorname{tw}}

\newcommand{\Ext}{{\operatorname{Ext}}}

\numberwithin{equation}{section}
\numberwithin{theo}{section}

\title{Lattice Point Asymptotics and Volume Growth on Teichm{\"u}ller Space.}
\author{Jayadev Athreya\thanks{partially supported by NSF grants DMS
    0603636, DMS 0244542 and DMS 1069153}, Alexander Bufetov\thanks{ 
A.I.B. is an Alfred P. Sloan Research Fellow.
He is supported in part by Grant MK-4893.2010.1 of the President of the Russian Federation, by the Programme on Mathematical Control Theory of the Presidium of the Russian Academy of Sciences,
by the Programme 2.1.1/5328 of the Russian Ministry of Education and Research,  by the Edgar Odell Lovett Fund at Rice University, by
the NSF under grant DMS~0604386, and by
the RFBR-CNRS grant 10-01-93115.
}, 
Alex Eskin\thanks{partially supported by NSF grants
    DMS 0244542, DMS 0604251 and DMS 0905912} 
and Maryam Mirzakhani\thanks{partially supported by the Clay
  foundation and by NSF grant DMS 0804136.}}

\begin{document}
\maketitle

\begin{abstract}
  We apply some of the ideas of the Ph.D. Thesis of G. A. Margulis
  \cite{margulis} to Teichm\"uller space. Let $X$ be a point in
  Teichm\"uller space, and let $B_R(X)$ be the ball of radius $R$
  centered at $X$ (with distances measured in the Teichm\"uller metric).
  We obtain asymptotic formulas as $R$ tends to infinity for the
  volume of $B_R(X)$, and also for the cardinality of the
  intersection of $B_R(X)$ with an orbit of the mapping class group.
\end{abstract}

\section{Introduction.}

Much of the study of the geometry and dynamics on Teichm\"uller and
moduli spaces is inspired by analogies with negatively curved spaces.
Two classical problems in the negative curvature setting are the
questions of \emph{lattice point counting} and \emph{volume growth
  entropy}. Let $M$ be a compact negatively curved Riemannian
manifold, and $\widetilde{M}$ its universal cover. $\Pi = \pi_1(M)$ acts
on $\widetilde{M}$ by isometries. Let $m_{BM}$ denote the Bowen-Margulis
measure on $\widetilde{M}$. Given $X, Y \in \widetilde{M}$, $R>0$, let
$B_R(X) \subset \widetilde{M}$ 
be the ball of radius $R$ centered at X. Let $p: \widetilde{M} \rightarrow M$ be
the natural projection. 

In his Ph.D. thesis
G.~A.~Margulis~\cite{margulis} showed:
\begin{theorem}\label{margthesis}  There is a function $c: M \times M
\rightarrow \bb R^+$ so that for every $X,Y \in
\widetilde{M}$, 

\begin{equation}\label{hyplattice}
|\Pi \cdot Y \cap B_R(X)| \sim c(p(X),p(Y)) e^{hR}
\end{equation}

\begin{equation}\label{hypvol}
m_{BM} (B_R(X)) \sim \left(\int_{M} c(p(X), z) dm_{BM}(z)\right) e^{hR}
\end{equation}

\noindent where $h> 0$ is the topological entropy of the geodesic
flow. Here and below, the notation $A \sim B$ means
that $A/B \to 1$ as $R \to \infty$.
\end{theorem}

In this paper, we apply some of the ideas of Margulis from~\cite{margulis} to
the problems of counting lattice points and understanding volume
growth entropy in Teichm{\"u}ller space. Our main results
(Theorems~\ref{theorem:asympt:count} and~\ref{theorem:volume:growth})
are analogues of Theorem~\ref{margthesis} for Teichm\"uller space, and
the fact that the Bowen-Margulis measure coincides with the Lebesgue
measure (see \S\ref{sec:conditional}). 
In the rest of this introduction we give the required
background and definitions of Teichm\"uller space and quadratic
differentials, and state our main results.

\subsection{Teichm\"uller space and quadratic differentials.}

We briefly recall some defintions. For a more detailed introduction,
see, e.g.,~\cite{FarbMargalit}. Let $g\geq 2$, and let $\Sigma_g$ be a
compact topological surface of genus $g$. Let $\cM_g$ and $\cT_g$
denote the moduli space and the Teichm{\"uller} space of compact
Riemann surfaces of genus $g$.

That is, $\cT_g$ is the space of equivalence classes of pairs $(X, f)$
where $X$ is a compact genus $g$ Riemann surface and $f: \Sigma_g
\rightarrow X$ is a diffeomorphism (known as a \emph{marking}). The
equivalence relation is given by $(X, f) \approx (Y, h)$ if there is a
biholomorphism $\phi: X \rightarrow Y$ so that $h^{-1} \circ \phi
\circ f$ is isotopic to the identity.

Let $\Gamma$ be the mapping class group of $\Sigma_g$, given by
isotopy classes of orientation preserving diffeomorphisms of
$\Sigma_g$. That is, $\Gamma =
\mbox{Diff}^+(\Sigma_g)/\mbox{Diff}_0^+(\Sigma_g)$. $\Gamma$ acts on
$\cT_g$ by changing the marking, and we have $\cM_g = \cT_g/\Gamma$.

$\cT_g$ carries a natural Finsler metric invariant under $\Gamma$, known as
the \emph{Teichm\"uller metric}. It is given by measuring the
quasi-conformal dilatation between surfaces, see,
e.g.~\cite{FarbMargalit}. The cotangent space of $\cT_{g}$ at a point
$X$ can be identified with the vector space $Q(X)$ of holomorphic
quadratic differentials on $X.$ Recall that given $X \in \cT_{g}$, a
quadratic differential $q \in Q(X)$ is a tensor locally given by
$\phi(z) dz^2$ where $\phi$ is holomorphic with respect to the
complex structure given by $X$.  Then the space
$\mathcal{Q}\te_{g}= \{(q,X) \;|\; X \in \te_{g}, q \in Q(X)\}$ is the
cotangent space of $\te_{g}.$ In this setting, the Teichm\"uller
metric corresponds to the norm
$$\parallel q \parallel_{\cT} = \int _{X} |\phi(z)|\; |dz|^{2}$$
on $\mathcal{Q}\te_{g}.$ 

Let $\qonemg$ (resp. $\qonetg$) denote the bundle of unit-norm
holomorphic quadratic differentials on $\cM_g$ (resp. $\cT_g$). We
have $\qonemg =\qonetg/\Gamma$.

The space $\qonetg$ carries a natural smooth measure $\mu$, 
preserved by the action of $\Gamma$   and 
such that $\mu(\qonemg) < \infty$ (\cite{masur}, \cite{veech}). 
We will fix a normalization for $\mu$ in \S\ref{sec:conditional}. 

We have a natural projection
$$
\pi: \qonetg \to \cT_{g}
$$
and we set $\mm=\pi_*\mu$.

\subsection{Statement of results.}

Let $X,Y\in \cT_g$ be arbitrary points, and 
let $B_R(X)$ be the ball of radius $R$ in the 
Teichm{\"u}ller space in the Teichm{\"u}ller metric, 
centered at the point $X$. 

Our main results are the following two theorems:

\begin{theorem}[Lattice Point Asymptotics]
\label{theorem:asympt:count}
As $R \to \infty$, 
\begin{displaymath}
|\Gamma \cdot Y \cap B_R(X)| \sim \frac{1}{h \, \mm(\cM_g)}\Lambda(X) \Lambda(Y) e^{h R},
\end{displaymath}
where $h = 6g-6$ is the entropy of the Teichm\"uller geodesic
flow with respect to Lebesgue measure (see \cite{veechteich}), 
and $\Lambda$ is a bounded function called the {\em
  Hubbard-Masur function}, which we define in
\S\ref{sec:hubbard:masur}.

\end{theorem}

\begin{theorem}[Volume Asymptotics]
\label{theorem:volume:growth}
As $R \to \infty$, 
$$
\mm(B_R(X)) \sim \frac{1}{h\, \mm(\cM_g)}e^{hR}\La(X)\cdot
\int_{\cM_g}\Lambda(Y)\, d\mm(Y).
$$
\end{theorem}

We also prove versions of the above theorems in ``sectors'', 
see Theorem~\ref{theorem:count:sectors}, 
Theorem~\ref{theorem:count:two:sectors}
and Theorem~\ref{theorem:volume:sectors} below.

\bold{Acknowledgments.}
The authors are grateful to Giovanni Forni, 
Vadim Kaimanovich, Yair Minsky  and Kasra Rafi for useful discussions, and especially to Howard Masur for his help
with all parts of the paper, in particular the appendix. 
We are also grateful to the Institute for Advanced Study, IHES 
and MSRI for their support. 

\subsection{Organization of the paper.}

This paper is organized as follows: In \S\ref{proofoutline}, we
describe the main ingredients in the proofs of our main theorems. We
construct foliations and measures satisfying the Margulis property of
uniform expansion in \S\ref{sec:conditional}. We define the
Hubbard-Masur function $\Lambda$ and describe its basic properties and
relation to counting multicurves in \S\ref{sec:hubbard:masur}. In
\S\ref{sec:measure:polar} we state a crucial lemma,
Proposition~\ref{prop:measure:polar},  on measure in polar
coordinates. We postpone the proof to \S\ref{sec:proof:measure:polar}.
In \S\ref{sec:mixing} we show how to use mixing to obtain
equidistribution results, and apply them to the counting problem in
\S\ref{sec:counting}. In \S\ref{sec:hodge:norm} we recall the
definition of the Hodge norm for abelian differentials
(\S\ref{sec:hodge:norm:abelian}), and extend and modify it to
quadratic differentials (\S\ref{sec:hodge:quadratic}). We use this to
define a distance along leaves of the stable and unstable foliations
for the Teichm\"uller geodesic flow, and compare it to other distances
on $\cT_g$ in \S\ref{sec:euclidean:hodge}. We prove a non-expansion
result for this distance,
Theorem~\ref{theorem:hodge:distance:nonincrease},  
in \S\ref{sec:hodge:distance}. In
\S\ref{sec:multiple:zeroes}, we use the results of
\S\ref{sec:hodge:norm} to prove the key estimates
Theorem~\ref{theorem:mz:volume} and
Theorem~\ref{theorem:badgeodesic:volume} 
from
\S\ref{sec:mixing}. Finally, in \S\ref{sec:volume}, we prove our
volume asymptotics result, Theorem~\ref{theorem:volume:growth},  
and relate it to counting multicurves.

\section{Outline of Proof.}\label{proofoutline}
\subsection{Notation and background.}
\label{sec:notation}
\bold{Teichm\"uller geodesic flow.}
We recall that when $g>1,$ the Teichm\"uller metric is not
Riemannian. However, geodesics in this metric are well understood.
A quadratic differential $q \in \mathcal{Q}\te_{g}$ with zeros at
$x_{1},\ldots, x_{k}$ is determined by an atlas of charts $\{
\phi_{i}\}$ mapping open subsets of $\Sigma_g-\{x_{1},\ldots,x_{k}\}$
to ${\Bbb R}^{2}$ such that the change of coordinates are of the form
$v \rightarrow \pm{v}+c.$ 
Therefore the group $\operatorname{SL}_{2}({\Bbb R})$ acts naturally on $\mathcal{Q}\cM_{g}$
by acting on the corresponding atlas; given $A \in
\operatorname{SL}_{2}({\Bbb R})$, $A \cdot q \in \mathcal{Q}\cM_{g}$
is determined by the new atlas $\{ A \phi_{i}\}.$

Let
$g_{t}=\begin{bmatrix} e^{t} &0\\
0 & e^{-t} \end{bmatrix}$. The action of the diagonal
subgroup $\{ g_t \,\mid\, t \in \reals \}$ 
is the Teichm\"uller geodesic flow for the Teichm\"uller metric (see
\cite{FarbMargalit}).

\bold{Warning.} In our normalization for the
  Teichm\"uller metric, the Teichm\"uller distance between $\pi(g_t
  q)$ and $\pi(q)$ is $t$. This normalization (and thus our value for
  the entropy $h$) differs by a factor of $2$ from that of e.g.
  \cite{veechteich}.  Our normalization is chosen in such a way that
  the top Lyapunov exponent of the {\it Kontsevich-Zorich cocycle} is
  equal to one.  In this case the top Lyapunov exponent of the flow is
  equal to {\it two}.  For a detailed discussion of the connection
  between the Lyapunov exponents of the Teichm\"uller geodesic flow
  and the Kontsevich-Zorich cocycle, see, e.g., \cite{forni}, where
  the same speed normalization is used for the Teichm\"uller geodesic
  flow as in our paper.  

\medskip

\noindent
We have \cite{veech}, \cite{masur}:
\begin{theo}{({\bf Veech, Masur})}
\label{theorem:VM}
The space $\cQ^{1}\cM_{g} $ carries a unique up to
  normalization measure $\mu$ 
in the Lebesgue measure class such that:
\begin{itemize} 
\item $\mu(\cQ^{1}\cM_{g}) < \infty$.
\item the action of $SL_{2}({\Bbb R})$ is volume preserving and ergodic;
\item the Teichm\"uller geodesic flow is mixing.
\end{itemize}
\end{theo}

\bold{Extremal lengths.}  Let $X$ be a Riemann surface. Then the
extremal length of a simple closed curve $\gamma$ on $X$ is defined by
$$
\label{eq:def:extremal:len}
\Ext_{\gamma}(X)=\sup_{\rho}
\frac{\ell_{\gamma}(\rho)^{2}}{\operatorname{Area}(X,\rho)},
$$ 
where the supremum is taken over all metrics $\rho$ conformally
equivalent to $X$, 
and $\ell_\gamma(\rho)$ denotes the length of $\gamma$ in the metric
$\rho$. 
The extremal length can be extended continuously from the space of
simple closed curves to the space of $\mf$ of measured foliations,
in such a way that $\Ext_{t\beta}(X) = t^2 \Ext_\beta(X)$ \cite{Kerckhoff:Asm}. 
On the other hand, by the uniformization theorem, each point $X \in \cT_{g}$ has a complete hyperbolic metric $\rho_{0}$ of constant 
curvature $-1$ in its conformal class. 
In general, for any simple closed curve $\alpha$,
 \begin{equation}
 \frac{\Ext_{\alpha}(X)} {\ell_{\alpha}(X)} \leq \frac{1}{2} e^{\ell_{\alpha}(X)/2},
\end{equation}
where $\ell_{\alpha}(X)$ is the length of the geodesic representative of $\alpha$ on $X$ with respect to the hyperbolic metric $\rho_{0}$ \cite{Maskit}.
Also, given $X$ there exists a constant $C_{X}$ such that 
$$\frac{1}{C_{X}} \ell_{\alpha}(X) \leq \sqrt{\Ext_{\alpha}(X)} \leq C_{X} \ell_{\alpha}(X). $$

The following result \cite{Kerckhoff:Asm} relates the ratios of extremal lengths to the Teichm\"uller distance:
\begin{theorem}({\bf Kerckhoff})
Given $X,Y \in \te_{g}$, the Teichm\"uller distance between $X$ and
$Y$ is given by 
$$d_{\cT}(X,Y)=\sup_{\beta \in \mathcal{C}}\log\left(\frac{\sqrt{\Ext_{\beta}(X)}}{\sqrt{\Ext_{\beta}(Y)}}\right),$$ 
where $\mathcal{C}$ is the set of simple closed curves on $\Sigma_g.$
\end{theorem} 
\subsection{The Margulis property.}
\label{sec:conditional}

Let $g_t$ be the Teichm{\"u}ller geodesic flow on $\qonetg$. Note that $g_t$
commutes with the action of $\Gamma$ and preserves the measure
$\mu$. We will be using the mixing property of the dynamical system
$(g_t, \qonemg, \mu)$ in \S\ref{sec:mixing}. 

Recall that a quadratic differential $q$ is uniquely determined by 
its imaginary and  real measured foliations $\eta^-(q)$ given by $\Im(q^{1/2})$
and $\eta^+(q)$ given by $\Re(q^{1/2})$. In this notation, we have $g_t q = g_t (\eta^+(q), \eta^-(q)) = (e^{t} \eta^+(q), e^{-t} \eta^-(q))$.
See e.g. \cite{FLP} for more details on measured foliations.

The flow $g_t$ preserves the following foliations:
\begin{enumerate}
\item $\F^{ss}$, whose leaves are sets of the form $\{q \st \eta^+(q)=const\}$;
\item $\F^{uu}$, whose leaves are sets of the form $\{q \st \eta^-(q)=const\}$.
\end{enumerate} 

In other words, for $q_0\in\qonetg$, a leaf of $\cF^{ss}$ is given by 
$$
\alpha^{ss}(q_0)=\{q\in\qonetg \st \ \eta^+(q)=\eta^+(q_0)\},
$$
and a leaf of $\cF^{uu}$ is given by
$$
\alpha^{uu}(q_0)=\{q\in\qonetg\st  \eta^-(q)=\eta^-(q_0)\}. 
$$
Note that the foliations $\F^{ss}$, $\F^{uu}$ are invariant under 
both $g_t$ and $\Gamma$; in particular, they 
descend to the moduli space $\qonemg$.

We also consider the foliations $\F^u$ whose leaves are defined by 
$$
\alpha^u(q)=\bigcup_{t\in{\mathbb R}} g_t\alpha^{uu}(q)
$$
and $\F^s$ whose leaves are defined by 
$$
\alpha^s(q)=\bigcup_{t\in{\mathbb R}} g_t\alpha^{ss}(q).
$$

Denote by $p:\mf\to\pmf$ the natural projection from the space 
of measured foliations onto the space of projective measured foliations. 
Now we may write 
$$
\alpha^{s}(q_0)=\{q\in\qonetg \st \ p(\eta^+(q))=p(\eta^+(q_0))\}; 
$$
$$
\alpha^{u}(q_0)=\{q\in\qonetg \st \ p(\eta^-(q))=p(\eta^-(q_0))\}. 
$$

The foliations $\F^u$ and $\F^{ss}$ form a complementary pair 
in the sense of Margulis 
\cite{margulis} (so do $\F^s$ and $\F^{uu}$, but the first pair will be more convenient for us).
Note that the foliations $\F^{ss}, \F^{s}, \F^u, \F^{uu}$ are, respectively, 
 the strongly stable, stable, unstable, 
and strongly unstable foliations for the Teichm{\"u}ller flow in the 
sense of Veech \cite{veechteich} and Forni \cite{forni}.
(See Theorem~\ref{theorem:hodge:distance:nonincrease} below for further
results in this direction). 

\bold{Conditional measures.}
The main observation, lying at the center of our construction, is the 
following. Each leaf $\alpha^{uu}$ of the 
foliation $\F^{uu}$, as well as each leaf $\alpha^{ss}$ of 
the foliation $\F^{ss}$ 
carries a globally defined normalized conditional measure $\mu_{\alpha^{ss}}$, invariant 
under the action the mapping class group, and having, moreover, the following property:
\begin{enumerate}
\item
  $(g_t)_*\mu_{\alpha^{uu}}=\exp(-ht)\mu_{g_t\alpha^{uu}}$; 
\item $(g_t)_*\mu_{\alpha^{ss}}=\exp(ht)\mu_{g_t\alpha^{ss}}$.
\end{enumerate}
where $h=6g-6$ is the entropy of the flow $g_t$ on $\qonemg$
with respect  to the smooth measure $\mu$.

The measures $\mu_{\alpha^{uu}}$ and $\mu_{\alpha^{ss}}$ may be
constructed as follows. 
Let $\nu$ denote the Thurston measure on $\mf$ \cite{FLP}. 
Note that each leaf $\alpha^s$ of $\cF^{s}$ is homeomorphic to an open
subset of 
$\mf$ via the map $\eta^-$. We can thus define the conditional measure on
this leaf to be the pullback of $\nu$, denote this measure by $\mu_{\alpha^s}$. Similarly one can define the
conditional measures $\mu_{\alpha^u}$ on leaves $\alpha^u$ of $\cF^{u}$. 

To define the measures on leaves of $\cF^{ss}$ (and $\cF^{uu}$), we restrict the conditional measure from a leaf  of
$\cF^{s}$ to a leaf of $\cF^{ss}$ (and similarly from a leaf of $\cF^{u}$
to a leaf of $\cF^{uu}$). This can be done explicitly in the following way:
for a subset $E \subset \mf$, let $Cone(E)$ denote the
cone based at the origin and ending at $E$ (i.e. the union of all the
line segments connecting the origin and points of $E$). We write
$\bar{\nu}(E)$ to denote $\nu(Cone(E))$. 
Now for a set $F \subset \alpha^{ss}$, 
we define $\mu_{\alpha^{ss}}(F) = \bar{\nu}(\eta^-(F))$. 
Similarly, for $F \subset \alpha^{uu}$, $\mu_{\alpha^{uu}}(F)$ is
defined to be $\bar{\nu}(\eta^+(F))$. 

In particular, if $\alpha_1$ and $\alpha_2$ are two leaves of the
foliation ${\cal F}^u$ and $U_1\subset \alpha_1$ and $U_2\subset
\alpha_2$ are chosen in such a way that $\eta^+(U_1)=\eta^+(U_2)$,
then we have $\mu_{\alpha_1}(U_1)=\mu_{\alpha_2}(U_2)$. The equality
$\eta^+(U_1)=\eta^+(U_2)$ is equivalent to the statement that $U_1$
may be taken to $U_2$ by holonomy along the leaves of the strongly
stable foliation ${\cal F}^{ss}$; the equality
$\mu_{\alpha_1}(U_1)=\mu_{\alpha_2}(U_2)$ thus means that the smooth
measure $\mu$ has the property of {\it holonomy invariance} with
respect to the pair of foliations $({\cal F}^u, {\cal F}^{ss})$.

This construction allows us to apply the arguments of G.A. Margulis 
and to compute the asymptotics for the volume of a ball of growing radius
in Teichm{\"u}ller space as well as the asymptotics of the number of
elements in the intersection of a ball with the orbit of
the mapping class group. The approach is similar, as noted above, to that of
\cite{EskMcM}.

\bold{Normalization of $\mu$.}
For convenience, we normalize the measure $\mu$ so that locally $d\mu =
d\mu_{\alpha^u}  d \mu_{\alpha^{ss}} = d\mu_{\alpha^s}
d\mu_{\alpha^{uu}}$.

\subsection{The Hubbard-Masur function.}
\label{sec:hubbard:masur}

\bold{The Hubbard-Masur Theorem.}
The Hubbard-Masur Theorem \cite{hubmas} states that given any point $X
\in \cT_g$ and any measured foliation $\beta \in \mf$, there exists 
a unique holomorphic quadratic differential $q$ on $X$ such that $\eta^+(q) =
\beta$. We also have the identity $\operatorname{Area}(q) = \Ext_{\beta}(X)$.

\bold{The measure $s_X$ and the multiple zero locus.}
For $X\in\cT_g$, we consider the unit (co)-tangent sphere 
$S(X)=\{q \in \qonetg \st \pi(q)=X\}$ at the point $X$. 
The conditional measure of $\mu$ on the sphere $S(X)$ will 
be denoted by $s_X$. It is by definition normalized so that $s_X(S(X))
= 1$. 

Let $\cP(1,\dots,1) \subset \qonetg$ denote the subset where the all
zeroes of the quadratic differential are distinct. $\cP(1,\dots,1)$ is
called the {\em principal stratum}. Its complement in $\qonetg$ is the
{\em multiple zero locus}. It is easy to see that the measures $s_X$
are defined for all $X$ and that this family of measures is smooth
away from the multiple zero locus.
We will also need the following:
\begin{theorem}
\label{theorem:sx:multiple:zero:locus}
For any $X \in \cT_g$, the measure $s_X$ gives zero weight to the
multiple zero locus. 
\end{theorem}

\bold{Proof.} See appendix A. 

\bold{The Hubbard-Masur function.}
Let $q\in\qonetg$ and let $\alpha^u(q)$ be the leaf of the foliation $\F^u$ containing $q$.
By the Hubbard-Masur Theorem, the projection 
$\pi$ induces a continuous bijection  between $\alpha^u(q)$ and
$\cT_g$ which is smooth away from the multiple zero locus. 
The mapping $\pi$ thus takes the
globally defined conditional measure $\mu_{\alpha^u(q)}$ on 
the fiber  $\alpha^u(q)$ to a measure on the 
Teichm{\"u}ller space; the resulting measure
$\pi_*(\mu_{\alpha^u(q)})$ is absolutely 
continuous with respect to the smooth measure $\mm$. Furthermore, by
Theorem~\ref{theorem:sx:multiple:zero:locus}, the measure $\mm$ is
also absolutely continuous with respect to $\pi_*(\mu_{\alpha^u(q)})$;
indeed, away from the multiple zero locus the mapping $\pi$ is smooth
with a smooth inverse, and the multiple zero locus has measure $0$ by
Theorem~\ref{theorem:sx:multiple:zero:locus}. 

We may therefore consider the corresponding Radon-Nikodym derivative.
Introduce a function $\lambda^+: \qonetg \to{\mathbb R}$ 
by the formula
$$
\lambda^+(q)=\frac{d\mm}{d(\pi_*(\mu_{\alpha^u(q)}))} \left( \pi(q)\right).
$$
Similarly, we define $\lambda^-: \qonetg\to{\mathbb R}$ via the formula
$$
\lambda^-(q)=\frac{d\mm}{d(\pi_*(\mu_{\alpha^s(q)}))}\left( \pi(q) \right).
$$

We set
$$
\Lambda(X)=\int_{S(X)} \lambda^+(q)\, ds_X(q) = \int_{S(X)} \lambda^-(q) \,
ds_X(q). 
$$
The equality of the two integrals will be justified by Proposition~\ref{prop:properties:Howard:Masur}. Note that the functions $\lambda^+, \lambda^-, \Lambda$ are $\Gamma$-invariant.
We call $\lambda^+$, $\lambda^-$ and $\Lambda$ the {\em Hubbard-Masur function}.

Note that by the Hubbard-Masur theorem,
$\eta^-$ (or $\eta^+$) defines a homeomorphism between the space of
all quadratic differentials at $X$  (with arbitrary area) and
the space $\mf$. This homeomorphism restricts to a homeomorphism
between $S(X)$ and  the set 
$$E_X = \{ \beta \in \mf \st \Ext_\beta(X) = 1  \}$$
where $\Ext_\beta(X)$ is the extremal length at $X$ of the measured
foliation $\beta$. Let $\delta^\pm_X: E_X \to S(X)$ denote the inverse of
$\eta^\pm$. 

It is easy to see that the functions
$\lambda^\pm$ are smooth on the 
complement of the multiple zero locus.

\begin{proposition}[Properties of $\lambda^+$, $\lambda^-$, $\Lambda$] 
\label{prop:properties:Howard:Masur} Let $X = \pi(q)$.
\begin{itemize}
\item[{\rm (i)}]
$\lambda^+(q) = \frac{d (\delta_X^-)_*(\bar{\nu})}{ds_X}(q)$, and $\lambda^-(q) = \frac{d (\delta_X^+)_*(\bar{\nu})}{ds_X}(q)$.
\item[{\rm (ii)}] $\Lambda(X) = \bar{\nu}(\eta^-(S(X))) = \bar{\nu}(\eta^+(S(X)))$.
\item[{\rm (iii)}] 
$\Lambda(X) = \bar{\nu}(E_X) = \nu(\{ \beta \in \mf \st \Ext_\beta(X) \le
1  \})$. 
\item[{\rm (iv)}] 
$\Lambda(X) = \bar{\nu}(\{ \beta \in \ml \st \Ext_\beta(X) \le
1  \})$, where $\ml$ is the space of measured laminations.
\end{itemize}
In (iv), by abuse of notation, $\nu$ denotes Thurston measure on $\ml$ and $\bar{\nu}(E)=\nu(Cone(E))$. 

\end{proposition}

\bold{Proof.} The property (i) follows from the fact that $d\mu =
d\mu_{\alpha^u}  d \mu_{\alpha^{ss}}$, the fact that $\eta^-_*(\alpha^{ss}) =
\bar{\nu}$, and 
the fact that if we write for $q \in \qonetg$, 
$q = (X,v)$ where $X = \pi(q)$ and $v \in S(X)$ then 
$d\mu(q) = d\mm(X) \, ds_X(v)$. \mc{MORE DETAILS}
The property (ii) follows from (i)
after making the change of variable $q \to \eta^-(q)$ in the definition
of $\Lambda$. The property (iii) follows from the fact that the image
$\eta^-(S(X))$ consists exactly of $E_X$. Finally (iv) is easily seen
to be equivalent to (iii). 
\qed\medskip

The following is proved in \S\ref{sec:volume}:
\begin{theorem}[Boundedness of $\Lambda$]
\label{theorem:Lambda:bounded}
There exists a constant $M$ such that 
$\Lambda(X) \le M $ for all $X \in \cT_g$.
\end{theorem}

For another interpretation of $\Lambda$ in terms of asymptotics of the
number of multicurves, see (\ref{eq:Lambda:and:multicurves}) below. 

\subsubsection{Conformal Densities}
Let the measure $\nu_X$  on $\pmf$ be the image of
$\lambda^-(q)ds_X(q)$ under the identification of $S(X)$ and
$\pmf$. 
That is, 
$$\nu_{X}(U)=\nu(\{\eta\; | \; [\eta]\in U\;, \Ext_{\eta}(X)\leq 1\}),$$\noindent where $\nu$ is, as above, Thurston measure on $\mf$, and $[\eta]$ denotes the image of $\eta$ in $\pmf$.
By definition, for $X,Y\in \te_{g}$, $\xi \in \mf$, we have 
$$
\frac{d\nu_X}{d\nu_Y}([\xi])=\left(\frac{\sqrt{{\rm Ext}_{\xi}(Y)}}{\sqrt{{\rm Ext}_{\xi}(X)}}\right)^{6g-6}.
$$
Note that the right-hand side only depends on $[\xi]$.
\begin{itemize}
\item
 Given $[\xi] \in \pmf$, let $\beta_{\xi}: \te_{g} \times \te_{g} \rightarrow \reals$ be the cocycle defined by 
$$\beta_{\xi}(X,Y)= \log\left(\frac{\sqrt{\Ext_{\xi}(X)}}{\sqrt{\Ext_{\xi}(Y)}}\right).$$
In our setting $\beta_{\xi}$ plays the role of the Busemann
cocycle for the Teichm\"uller metric, and (formally) the family
$\{\nu_{X}\}_{X \in \te_{g}} $ of finite measures on $P\MF$ is a
family of {\it conformal densities} of dimension $\delta=6g-6$ for the cocycle
$\beta$; i.e, 
for any $X,Y \in \te_{g}$ and almost every $\xi \in \pmf$ we have
$$\frac{d\nu_{X}}{d\nu_{Y}}([\xi])=\operatorname{exp}(\delta \beta_{\xi}(Y,X)), \;\;
\mbox{and} \;\; \nu_{\gamma \cdot X}(\gamma \cdot U)= \nu_{X}(U),$$
where $U \subset \pmf$, and $\gamma$ is an element of the mapping class group.

Moreover, by the measure classification result of \cite{LM},
$\{\nu_{X}\}_{X \in \te_{g}}$ is the unique conformal density for the
action of the mapping class group on $\te_{g}$ up to scale;  
the only possibility for $\delta$ is $6g-6.$

\item Given $X \in \te_g$, and transverse minimal foliations $\xi,
\eta\in \mf$, define 
$$
\beta(X, [\xi], [\eta])=\left(\frac{\sqrt{{\rm Ext}_{\xi}(X){\rm Ext}_{\eta}(X)}}{A(\xi, \eta)}\right)^{6g-6},
$$
 where $A(\xi, \eta)$ is the area of the quadratic 
differential with horizontal foliation $\xi$, and vertical foliation $\eta$. Note that the right-hand side only depends on $[\xi]$, $[\eta]$.

Now introduce a measure ${\mathbb P}$ on the product $\pmf\times \pmf$
by the formula
\begin{equation}
{\mathbb P}=\beta(X, [\xi], [\eta])\nu_X([\xi])\nu_X([\eta]).
\end{equation}
Observe that the right-hand side does not depend on $X$ and is invariant 
under the diagonal action of the mapping class group.

The space $\qonetg$ of quadratic differentials fibres over $\pmf\times
\pmf$ (with fibre ${\mathbb R}$) and the measure ${\mathbb P}$ on
$\pmf\times \pmf$ lifts to (a multiple of) the Lebesgue measure on
$\qonetg$ (by construction, it is invariant and belongs to the Lebesgue
measure class --- but this determines the measure uniquely).
\end{itemize}

\noindent A consequence of Theorem~\ref{theorem:count:sectors} is that the measures $\nu_X$ are the \emph{Patterson-Sullivan} measures on $\cT_g$. These constructions are similar to the ones due to Kaimanovich in \cite{kaim1} and \cite{kaim2}.

\subsection{Measure in Polar Coordinates.}
\label{sec:measure:polar}


Recall that $\qonetg$ is the space of pairs $(X,q)$ where $X$ is
a Riemann surface and $q$ is a holomorphic quadratic differential on
$X$. (We occasionally drop $X$ from the notation and denote points
of $\qonetg$ by $q$ alone). 


Fix $X \in \cT_g$. Then we have a diffeomorphism $\Phi: S(X) \cross
\reals^+ \to \qonetg$ where $\Phi(q,t) = g_t q$ is the point in $\qonetg$
which is the endpoint of the length $t$ geodesic segment starting at
$X$ and tangent to $q \in S(X)$. We can then write our measure $d\mm$ in polar coordinates,
\begin{equation}
\label{eq:def:Delta}
d\mm(\pi(g_t q))  = \Delta(q,t) \, ds_X(q) \, dt.
\end{equation}

\medskip

\noindent\textbf{Remark:} To be completely precise, we should write $d\mm(\pi(g_t q))  = \Phi_*( \Delta(q,t) \, ds_X(q) \, dt)$, but it is standard practice to avoid this formality with polar coordinates.

\medskip

\begin{proposition}
\label{prop:measure:polar}
Let $K_1 \supset K$ be compact subsets of $\cP(1,\dots,1)$, $q \in 
\qonetg$, and
suppose $q$ and $g_t q$ both belong to $\Gamma K$. 
Then 
\begin{equation}
\label{eq:measure:polar:crude:bound}
|\Delta(q,t)| \le C e^{ht}, 
\end{equation}
where $C$ depends only on $K$. If in addition 
\begin{equation}
\label{eq:s:half:time}
|\{ s \in [0,t] \st g_s q \in \Gamma K_1 \} | \ge (1/2) t
\end{equation}
then there exists $\alpha > 0$ so that
\begin{equation}
\label{eq:measure:polar:asymptiotics}
\Delta(q,t) = e^{ht} \lambda^-(q) \lambda^+(g_t q) + O(e^{(h-\alpha)t})
\end{equation}
where $\alpha$ and the implied constant on
(\ref{eq:measure:polar:asymptiotics}) depend only on $K_1$. \mc{need
  $K1$?} \mc{check $\lambda^+$ vs $\lambda^-$?}
\end{proposition}

\medskip

\noindent We will prove this proposition in \S\ref{sec:proof:measure:polar}.

\subsection{Mixing.}
\label{sec:mixing}

Let $\cU$ be a open subset of the boundary at infinity of $\cT_g$
(which can be identified with the space $\cP\cM\cF$ of projective
measured foliations). For
each $X \in \cT_g$ we may identify  $\cU$ with a subset $\cU(X)$ of $S(X)$ \mc{
  beware non-uniquely-ergodic}. Let $Sect_\cU(X)
\subset \cT_g$ denote the set $\bigcup_{t \ge 0} \pi(g_t \cU(X))$, i.e. the
sector based at $X$ in the direction $\cU$. 
For a subset $A \subset \cT_g$, 
$Nbhd_{r}(A)$ denotes the set of points within Teichm\"uller
distance $r$ of $A$.

\begin{theorem}
\label{theorem:mz:volume}
For any $X \in \cT_g$ and any $\epsilon > 0$ there exists an open
subset $\cU$ of $S(X)$ containing the intersection of $S(X)$ with the 
multiple zero locus such that for any compact $\cK \subset \cT_g$
\begin{equation}
\label{eq:mz:volume}
\limsup_{R \to \infty} \ e^{-hR} \mm( Nbhd_1(B_R(X) \cap Sect_\cU(X)
\cap \Gamma\cK))  \le \epsilon.
\end{equation}
\end{theorem}

Let $\cK$ be a compact subset of $\cT_g$, and let $K'$ be a compact
subset of the principal stratum $\cP(1,\dots,1) \subset \qonetg$. Let
$B_R(X,\cK,K')$ denote the set of points $Y \in \Gamma \cK$ such that
$d(X,Y) < R$, and the geodesic from $X$ to $Y$ spends at least half
the time outside $\Gamma K'$.

\begin{theorem}
\label{theorem:badgeodesic:volume}
For $\epsilon > 0$ and any compact $\cK \subset \cT_g$ there exists
compact $K' \subset \cP(1,\dots,1)\subset \qonetg$ such that for any $X
\in \cK$,
\begin{equation}
\label{eq:badgeodesic:volume}
\limsup_{R \to \infty} \ e^{-hR} \mm( Nbhd_1(B_R(X,\cK, K')))
\le \epsilon.
\end{equation}
\end{theorem}

Theorems~\ref{theorem:mz:volume} and ~\ref{theorem:badgeodesic:volume} will be proved in \S\ref{sec:multiple:zeroes}.

Let $\psi$ be a nonnegative compactly supported continuous function on
$\cM_g$, which we extend to a function $\psi: \qonemg \to
\reals$ by 
making it constant on each sphere $S(X)$. We can consider $\psi$ to be
a $\Gamma$-periodic function on $\cT_g$ (or $\qonetg$). Let $\phi$ be
another such function.

\begin{proposition}[Mean Ergodic Theorem]
\label{prop:mean:ergodic:theorem}
We have 
\begin{eqnarray}
\label{eq:mean:ergodic}
\lim_{R \to \infty} h e^{-hR} \int_{\cM_g} \phi(X) \left( \int_{B_R(X)}
  \psi(Y) \, 
  d\mm(Y) \right) \, d\mm(X)  \\ = \frac{1}{\mm(\cM_g)} \int_{\cM_g} \phi(X) \Lambda(X) \,
d\mm(X) \int_{\cM_g} \psi(Y) \Lambda(Y) \, d\mm(Y).\nonumber
\end{eqnarray}
\end{proposition}

\bold{Proof.} 
Suppose $\epsilon > 0$ is given.
Let $\cK$ be the union
of the supports of $\phi$ and $\psi$. 
Without loss of generality, we may assume that the support of $\phi$
is small enough so that there exists an open $\cU \subset \partial
\cT_g$ such that for each $X$ in the support of $\phi$, $\cU(X)$
contains the intersection of $S(X)$ with the multiple zero locus, and
also (\ref{eq:mz:volume}) holds. Here, as above, $\cU(X)$ is the sector in $S(X)$ corresponding to directions in $\cU$.

We break up the integral over $B_R(X)$ in (\ref{eq:mean:ergodic}) into
two parts: the first over $B_R(X) \cap Sect_{\cU}(X) \cap \Gamma\cK$
and the second over the complement. The first is bounded by $C_1
\epsilon$ in view of Theorem~\ref{theorem:badgeodesic:volume}.

Let $K'$ be as in Theorem~\ref{theorem:badgeodesic:volume}.
We use polar coordinates for the  integral over
$B_R(X)$. We get, 
\begin{eqnarray}
\label{eq:use:polar}
\int_{\cM_g} \phi(X) \left( \int_{B_R(X)}
  \psi(Y) \, 
  d\mm(Y) \right) \, d\mm(X) =\\ \int_{\cM_g} \phi(X) \left( \int_0^R
  \int_{S(X)} 
  \Delta(q,t) \psi(g_t q) \, ds_X(q) \, dt\right) \, d\mm(X).\nonumber
\end{eqnarray}  

In view of (\ref{eq:badgeodesic:volume}), we may assume, (up to error
of size $C \epsilon e^{hR}$) that $q \in S(X)$ belongs to the compact
set $\Gamma K'' = \bigcup_{X \in \cK} S(X) \setminus \cU(X)$, which is
away from 
the multiple zero locus. 
Also note that the left-hand-side of (\ref{eq:use:polar}) is
symmetric in $X$ and $Y$ (up to interchanging the functions $\phi$ and
$\psi$). So we may also assume that $g_t q$ also belongs to $\Gamma K''$. 
Also in view
of (\ref{eq:badgeodesic:volume}), we may assume that
(\ref{eq:s:half:time}) holds. (Again the contribution of the part of the
integral where this assumption is violated is bounded by $C \epsilon e^{hR}$.) 

We now consider the part of the region where none of the three assumptions are
violated, i.e. $q \in \Gamma K''$, $g_t q \in \Gamma K''$ and
(\ref{eq:s:half:time}) holds. 
We are finally in a position to use 
Proposition~\ref{prop:measure:polar} to replace $\Delta(q,t)$ by 
$e^{ht}\lambda^-(q) \lambda^+(g_t q)$. Let $\widehat{\lambda^{\pm}} $
denote the truncations of $\lambda^{\pm}$ to $\Gamma K''$. We can write
\begin{multline}
\label{eq:triple:int}
\int_0^R e^{ht}\int_{\cM_g} \phi(X) \int_{S(X)} \widehat{\lambda^-}(q)
\psi(g_t q) \widehat{\lambda^+}(g_t q) \, ds_X(q) \, d\mm(X) \, dt  =
\\ 
\int_0^R  e^{ht} \int_{\qonemg} \phi(q) \widehat{ \lambda^-}(q)
\psi(g_t q) \widehat{\lambda^+}(g_t q)  \, d\mu(q) \, dt, 
\end{multline}
where we used the fact that $\phi$ is constant on each sphere $S(X)$
and the formula $$d\mu(q) = d \mm(\pi(q)) \, d{s_{\pi(q)}}(q).$$ 
Recall that $\mu(\qonemg) < \infty$ by Theorem~\ref{theorem:VM}. 
Now we can apply the mixing property of the geodesic flow to the
functions $\phi \widehat{\lambda^-}$ and $\psi \widehat{ \lambda^+}$. The right-hand side
of (\ref{eq:triple:int}) is then asymptotic to 
\begin{multline*}
\frac{e^{hR}}{h\, \mu(\qonemg)} \int_{\qonemg} \phi(q)\widehat{ \lambda^-}(q) \, d\mu(q) \int_{\qonemg}
\psi(q) \widehat{ \lambda^+}(q) \, d\mu(q)  \\ = \frac{e^{hR}}{h\, \mm(\cM_g)}
\int_{\cM_g} \phi(X) \Lambda(X) \, d\mm(X)
\int_{\cM_g} \psi(Y) \Lambda(Y) \, d\mm(Y), 
\end{multline*}
where the equality holds up to the truncation error. So the right hand side of (\ref{eq:mean:ergodic}) is equal to the left
hand side of (\ref{eq:mean:ergodic}) up to an error which is bounded by $C \epsilon$, where $C$ depends only on $\phi$ and $\psi$. Since $\epsilon > 0$ is
arbitrary, the theorem follows. 
\qed\medskip

The proof of Proposition~\ref{prop:mean:ergodic:theorem} also shows
that for any open $\cU \subset S(X)$, 
\begin{multline}
\label{eq:mean:sectors}
\lim_{R \to \infty} h e^{-hR} \int_{\cM_g} \phi(X) \left(
  \int_{B_R(X) \cap Sect_{\cU}(X)} \psi(Y) \, 
  d\mm(Y) \right) \, d\mm(X)  = \\ 
\frac{1}{\mm(\cM_g)}\int_{\cM_g} \phi(X) \left(\int_{\cU}
\lambda^-(q) \, ds_X(q) \right) \,
d\mm(X) \int_{\cM_g} \psi(Y) \Lambda(Y) \, d\mm(Y).
\end{multline}

\subsection{Counting.}
\label{sec:counting}
\bold{Proof of Theorem~\ref{theorem:asympt:count}.}
Let $F_R(X,Y)$ denote the cardinality of the intersection of $B_R(X)$
with $\Gamma \cdot Y$. Let $\phi$ and $\psi$ be non-negative
compactly supported continuous functions on
$\cM_g$. Note that
\begin{displaymath}
\int_{\cM_g} F_R(X,Y) \psi(Y) \, d\mm(Y) = \int_{B_R(X)}
\psi(Y) \, d\mm(Y).
\end{displaymath}
Hence, by Proposition~\ref{prop:mean:ergodic:theorem}, 
\begin{multline*}
\lim_{R \to \infty} e^{-hR} \int_{\cM_g} \phi(X) \int_{\cM_g} F_R(X,Y)
\psi(Y) \, d\mm(Y) \, d\mm(X) \\ =  
\frac{1}{h\, \mm(\cM_g)} \int_{\cM_g}\phi(X) \Lambda(X) \, d\mm(X) 
\int_{\cM_g} \psi(Y) \Lambda(Y) \, d\mm(Y), 
\end{multline*}
i.e. $e^{-hR} F_R(X,Y)$ converges weakly to $\frac{1}{h \mm(\cM_g)}\Lambda(X)
\Lambda(Y)$. 
We would like to show that the convergence is pointwise and uniform on compact
sets.

Let $\epsilon > 0$ be arbitrary. 
Let $\psi$ be  supported on a ball of radius $\epsilon$ around $Y$ (in
the Teichm\"uller distance),
and satisfy $\int_{\cM_g} \psi\, d\mm = 1$. Let $\phi$ be
supported on a ball of radius  $\epsilon$ around $X$, with
$\int_{\cM_g} \phi\, d\mm = 1$. 
By the triangle inequality, \mc{GIVE MORE DETAILS}
\begin{displaymath}
\int_{\cM_g} \phi(X) \int_{B_{R-3\epsilon}(X)}
\psi(Y) \, d\mm(Y) \, d\mm(X) \le F_R(X,Y) \le \int_{\cM_g}
\phi(X) \int_{B_{R+3\epsilon}(X)} \psi(Y) \, d\mm(Y) \, d\mm(X). 
\end{displaymath}
After multiplying both sides by $h e^{-hR}$ and taking $R \to \infty$,
we get, after applying Proposition~\ref{prop:mean:ergodic:theorem},
\begin{multline*}
e^{-3 h \epsilon} \frac{1}{\mm(\cM_g)} \int_{\cM_g} \phi(X) \Lambda(X)
\, d\mm(X) \int_{\cM_g} 
\psi(Y) \Lambda(Y) \, d\mm(Y) 
\le \liminf_{R \to \infty} h e^{-h R} F_R(X,Y) \le \\ \le \limsup_{R \to
  \infty} h e^{-h R} F_R(X,Y) \le e^{3 h \epsilon}
\frac{1}{\mm(\cM_g)}\int_{\cM_g} 
\phi(X) \Lambda(X) \, d\mm(X) \int_{\cM_g} \psi(Y) \Lambda(Y)
\, d\mm(Y). 
\end{multline*}
Since $\epsilon > 0$ is arbitrary, and $\Lambda$ is uniformly
continuous on compact sets, we get
\begin{displaymath}
\lim_{R \to \infty} h e^{-h R} F_R(X,Y) = \frac{1}{\mm(\cM_g)}
\Lambda(X) \Lambda(Y) 
\end{displaymath}
as required. 
\qed\medskip

In fact, the same argument using (\ref{eq:mean:sectors}) instead of
Proposition~\ref{prop:mean:ergodic:theorem}
 yields the following:
\begin{theorem}
\label{theorem:count:sectors}
Suppose $X \in \cT_g$, and $\cU$ is an open subset of the boundary
at infinity of $\cT_g$. Then for any $Y \in
\cT_g$, we have, as $R \to \infty$, 
\begin{displaymath}
|B_R(X) \cap Sect_\cU(X) \cap \Gamma \cdot Y| \sim \frac{1}{h\,
  \mm(\cM_g)}e^{h R} \Lambda(Y) \int_\cU 
\lambda^-(q) \, ds_X(q).
\end{displaymath}
\end{theorem}

We also obtain the following:
\begin{theorem}
\label{theorem:count:two:sectors}
Suppose $X \in \cT_g$, and $\cU$ and $\cV$ are open subsets of
the boundary at infinity of $\cT_g$. Then for any $Y \in
\cT_g$, we have, as $R \to \infty$, 
\begin{multline*}
|\{ \gamma \in \Gamma \st \gamma\; \cdot Y \in B_R(X) \cap Sect_\cU(X) \text{
  and } \gamma^{-1} \cdot  X \in Sect_{\cV}(Y) \} | \\ 
\sim \frac{1}{h\,
  \mm(\cM_g) }e^{h R}
\int_\cU \lambda^-(q) \, ds_X(q) \int_\cV \lambda^+(q) \, ds_Y(q).
\end{multline*}
\end{theorem}

The proof of Theorem~\ref{theorem:count:two:sectors} is very similar
to that of Theorem~\ref{theorem:count:sectors} except that we do
not assume that the function $\psi: \qonemg \to \reals$ 
is the pullback under $\pi$ of a function from $\cM_g \to
\reals$. The details are left to the reader.


\section{The Hodge Norm.}
\label{sec:hodge:norm}

\subsection{The Hodge Norm for Abelian Differentials.}
\label{sec:hodge:norm:abelian}

\subsubsection{Definition and Basic Properties.}
\label{sec:hodge:def}
Let $X$ be a Riemann surface. By definition, $X$ has a complex
structure. Let $\cH_X$ denote the set of holomorphic $1$-forms on
$X$. One can define {\em Hodge inner product} on $\cH_X$ by
\begin{displaymath}
\langle \omega, \eta \rangle = \int_X \omega \wedge \bar{\eta}.
\end{displaymath}
We have a natural map $r: H^1(X,\reals) \to \cH_X$ which sends
a cohomology class $\lambda \in H^1(X,\reals)$ to the 
holomorphic $1$-form $r(\lambda) \in \cH_X$ such that 
the real part of $r(\lambda)$ (which is a harmonic $1$-form)
represents $\lambda$. We can thus define  the Hodge inner product on
$H^1(X,\reals)$ by $\langle \lambda_1, \lambda_2 \rangle = \langle
r(\lambda_1), r(\lambda_2) \rangle$. We have
\begin{displaymath}
\langle \lambda_1, \lambda_2 \rangle = \int_X \lambda_1 \wedge *\lambda_2,
\end{displaymath}
where $*$ denotes the Hodge star operator, and we choose harmonic representatives of $\lambda_1$ and $*\lambda_2$ to evaluate the integral.
We denote the associated norm by $\| \cdot \|$. This is the {\em Hodge
  norm}. See \cite{FarkasKra}.

Let $\alpha$ be a homology class in $H_1(X,\reals)$. We can define the
cohomology class $*c_{\alpha} \in H^1(X,\reals)$ so that for all
$\omega \in H^1(X, \reals)$, 
\begin{displaymath}
\int_\alpha \omega = \int_X \omega \wedge *c_\alpha. 
\end{displaymath}
Then, 
\begin{displaymath}
\int_X *c_\alpha \wedge *c_\beta = I(\alpha,\beta),
\end{displaymath}
where $I(\cdot, \cdot)$ denotes algebraic intersection number.
We have, for any $\omega \in H^1(X,\reals)$, 
\begin{displaymath}
\langle \omega, c_\alpha \rangle = \int_X \omega \wedge *c_\alpha = 
\int_\alpha \omega.
\end{displaymath}
We note that $*c_\alpha$ is a purely topological construction which 
depends only on $\alpha$, but $c_\alpha$
depends also on the complex structure of $X$.

\subsubsection{The Hodge Norm and the Hyperbolic Metric.}
\label{sec:hodge:hyp}
Let $\alpha$ be a simple close curve on a Riemann surface $X$. Let
$\ell_\alpha(\sigma)$ denote the length of the geodesic 
representative of $\alpha$ in the hyperbolic metric which is in the
conformal class of $X$. 

Fix $\epsilon_* > 0$ (the \emph{Margulis constant}) so that any two
curves of hyperbolic length less than $\epsilon_*$ must be disjoint. 

\begin{theorem}
\label{theorem:hodge:hyperbolic}
For any constant $D > 1$ there exists a constant $c > 1$ such that
for any simple closed curve
$\alpha$ with $\ell_\alpha(\sigma) < D$, 
\begin{equation}
\label{eq:hodge:hyperbolic}
\frac{1}{c} \ell_\alpha(\sigma)^{1/2} \le \| c_\alpha \| <
c \, \ell_\alpha(\sigma)^{1/2}.
\end{equation}
Furthermore, if $\ell_\alpha(\sigma) < \epsilon_0$ and $\beta$ is the
shortest simple closed curve crossing $\alpha$, then
\begin{displaymath}
\frac{1}{c} \ell_\alpha(\sigma)^{-1/2} \le \| c_\beta \| <
c \, \ell_\alpha(\sigma)^{-1/2}.
\end{displaymath}
\end{theorem}

\bold{Proof.}  Let $\alpha_1, \dots, \alpha_n$ be the curves with
hyperbolic length less than $\epsilon_0$. For $1 \le k \le n$, let
$\beta_k$ be the shortest curve with $i(\alpha_k, \beta_k) =1$, where
$i(\cdot,\cdot)$ denotes the geometric intersection number. It is
enough to prove (\ref{eq:hodge:hyperbolic}) for the $\alpha_k$ and the
$\beta_k$ (the estimate for other moderate length curves follows from
a compactness argument).
 
We can find a collar region around $\alpha_k$ as follows:
take two annuli $\{ z_k \st 1 > | z_k | > |t_k|^{1/2}\}$ and $\{ w_k \st 1 > w_k> |t_k|^{1/2} \}$ and identify the inner boundaries via the map $ w_k =
t_k/z_k$. (This coordinate system is used in e.g. \cite[Chapter 3]{fay},
also \cite{Masur:WP}, \cite{forni}, \cite[\S{3}]{wolpert} and elsewhere). 
The hyperbolic metric $\sigma$
in the collar region is approximately $|dz|/(|z| |\log |z||)$. 
Then $\ell_{\alpha_k}(\sigma) \approx 1/|\log t_k|$. 
By \cite[Chapter 3]{fay}
any holomorphic $1$-form $\omega$ can be written in
the collar region as
\begin{displaymath}
\left( a_0(z_k+t_k/z_k,t_k) + \frac{a_1(z_k+t_k/z_k,t_k)}{z_k} \right)
\, dz_k, 
\end{displaymath}
where $a_0$ and $a_1$ are holomorphic in both variables. (We assume
here that the limit surface on the boundary of Teichm\"uller space is
fixed). This implies that as $t_k \to 0$, 
\begin{displaymath}
\omega = \left(\frac{a}{z_k} + h(z_k)+ O(t_k/z_k^2)\right) \, dz_k
\end{displaymath}
where $h$ is a holomorphic function which 
remains bounded as $t_k \to 0$, and the implied constant is bounded as
$t_k \to 0$. (Note that when $|z_k| \ge |t_k|^{1/2}$, $|t_k/z_k^2| \le 1$). 
Now from the condition
$\int_{\alpha_k} *c_{\beta_k} = 1$ we see that on the collar of
$\alpha_j$, 
\begin{equation}
\label{eq:c:beta:holomorphic}
c_{\beta_k} + i*\!\!c_{\beta_k} = \left(\frac{\delta_{kj}}{(2 \pi) z_j} +
h_{kj}(z_j) + O(t_j/z_j^2)\right) \, dz_j, 
\end{equation}
where the $h_{kj}$  are  holomorphic and bounded as $t_j \to 0$. (We use the notation $\delta_{kj} = 1$ if $k=j$ and zero otherwise). 
Also from the condition $\int_{\beta_k} *c_{\alpha_k} = 1$ we have
\begin{equation}
\label{eq:c:alpha:holomorphic}
c_{\alpha_k} + i*\!\!c_{\alpha_k} = \frac{1}{|\log t_j|}\left(\frac{\delta_{kj}}{z_j} +
s_{kj}(z_j) + O(t_j/z_k^2)\right) \, dz_j, 
\end{equation}
where $s_{kj}$ also remains holomorphic and is bounded as $t_j \to 0$. 
Then,
\begin{displaymath}
\|c_{\beta_k}\|^2 = \|*\!\!c_{\beta_k}\|^2 = \int_X c_{\beta_k} \wedge *
c_{\beta_k} \approx \int_0^{2\pi} \int_{t_k}^1 \frac{1}{4 \pi^2 r^2} r
\, dr \, d\theta \approx \frac{|\log t_k|}{4 \pi}
\end{displaymath}
and
\begin{displaymath}
\|c_{\alpha_k}\|^2 = \|*\!\!c_{\alpha_k}\|^2 = \int_X c_{\alpha_k} \wedge *
c_{\alpha_k} \approx \int_0^{2\pi} \int_{t_k}^1 \frac{1}{(\log t_k)^2 r^2} r
\, dr \, d\theta \approx \frac{2 \pi}{|\log t_k|}.
\end{displaymath}
\qed\medskip

\subsubsection{The Hodge Norm and Extremal Length.}\label{sec:hodge:norm:extremal}
We recall the following theorem (\cite{accola}, \cite{blatter}):
\begin{theorem}[Accola, Blatter]
\label{theorem:accola:blatter}
For any Riemann surface $X$ and any $\lambda \in H_1(X,\reals)$,
\begin{displaymath}
\|c_\lambda\|^2 = \sup_{\rho} \inf_{\gamma \in [\lambda]} \frac{\ell_{\gamma}(\rho)^2}{\operatorname{Area}(\rho)},
\end{displaymath}
where the supremum is over the metrics $\rho$ which are in the
conformal class of $X$, and the infimum is over all the multicurves
 $\gamma$ homologous to $\lambda$. 
\end{theorem}
Note that in the definition (\ref{eq:def:extremal:len}) of extremal
length, the infimum is 
is over all the simple closed curves $\gamma$ {\em homotopic} to $\lambda$. Recall
that the extremal metric $\rho$ 
is always the flat metric obtained from a quadratic differential $q$,
and in this metric the entire surface consists of a flat cylinder in
which $\lambda$ is the core curve. \mc{add reference}
Thus, if this quadratic differential $q$ 
is in fact the square of an Abelian differential
then it follows from Theorem~\ref{theorem:accola:blatter}
that
\begin{equation}
\label{eq:extremal:hodge}
\Ext_{\lambda}(X) = \|c_\lambda\|^2.
\end{equation}

\subsubsection{The Hodge Norm and the Geodesic Flow.}
\label{sec:hodge:geod}
Let $\omegaonetg$ denote the space of (unit area) holomorphic $1$-forms on surfaces of
genus $g$. Recall that $g_t$, the Teichm\"uller geodesic flow, preserves $\omegaonetg$
(where we map $\omegaonetg$ into $\qonetg$ by squaring abelian differentials). Fix a point $S$ in $\omegaonetg$; then $S$ is a pair
$(X,\omega)$ where $\omega$ is a holomorphic $1$-form on $X$. Let $\|
\cdot \|_0$ denotes the Hodge norm on the surface $X_0 = X$, and let
$\| \cdot \|_t$ denote the Hodge norm on the surface $X_t = \pi(g_t
S)$.

The following fundamental result is due to Forni \cite[\S{2}]{forni}: 
\begin{theorem}
\label{theorem:forni}

For any $\theta \in H^1(X,\reals)$ and any $t \ge 0$, 
\begin{displaymath}
\|\theta\|_t \le  e^{t} \| \theta\|_0.
\end{displaymath}
Moreover, there exists a constant $C$ depending only on the genus, such that 
if  $\langle \theta,\omega\rangle =0 $, and 
for some compact subset $\cK$ of $\cM_g$, 
the geodesic segment $[S,g_t S]$ spends at least half the
time in $\pi^{-1}(\cK)$, then we have
\begin{displaymath}
\|\theta\|_t \le C e^{(1-\alpha)t} \| \theta\|_0,
\end{displaymath}
where $\alpha > 0$ depends only on $\cK$.
\end{theorem}

\subsection{Proof of Proposition~\ref{prop:measure:polar}.}
\label{sec:proof:measure:polar}

In this section we prove Proposition~\ref{prop:measure:polar}. We recall the setting: given $X \in \cT_g$, let $g_t q$ be the point in $\qonetg$
which is the endpoint of the length $t$ geodesic segment starting at
$X$ and tangent to $q \in S(X)$. $\Delta(q, t)$ is then given by:

$$d\mm(\pi(g_t q))  = \Delta(q,t) \, ds_X(q) \, dt.$$

\noindent Let $K_1 \supset K$ be compact subsets of $\cP(1,\dots,1)$, $q \in 
\qonetg$, and suppose $q$ and $g_t q$ both belong to $\Gamma K$. Proposition~\ref{prop:measure:polar} asserts a general bound, inequality (\ref{eq:measure:polar:crude:bound}):

$$|\Delta(q,t)| \le C e^{ht}$$

\noindent where $C$ depends only on $K$. If we are given some more information about the trajectory $\{g_s q\}_{0 \le s \le t}$, in particular that it spends a portion of its time in a compact set, we can say more. Precisely, if the inequality (\ref{eq:s:half:time}):

$$|\{ s \in [0,t] \st g_s q \in \Gamma K_1 \} | \ge (1/2) t$$

\noindent is satisfied then there exists $\alpha > 0$ so that the inequality (\ref{eq:measure:polar:asymptiotics}):

$$\Delta(q,t) = e^{ht} \lambda^-(q) \lambda^+(g_t q) + O(e^{(h-\alpha)t})$$

\noindent holds. Moreover, $\alpha$ and the implied constant in
(\ref{eq:measure:polar:asymptiotics}) depend only on $K_1$.

\bold{Proof.} Without loss of generality, $q \in K$.  Let $X =
\pi(q)$, $Y = \pi(g_t q)$.  Write $Y = \gamma z$ where $\gamma \in
\Gamma$ and $z \in K$.  Let $A$ be the matrix of $\gamma$ acting on
the odd part of the homology of the orienting double cover $\tilde{X}$
of $q$ (i.e.  $A$ is the Kontsevitch-Zorich cocycle.)  Let $\langle
\cdot, \cdot \rangle$ denote the Hodge inner product, and $A^*$ denote
the adjoint of $A$, where we view $A$ as a map between the inner
product spaces determined by the Hodge inner product at $X$ and $Y$
respectively. Let $n = h = \dim H^1_{odd}(\tilde{X},\reals)$, and let
$e_1, \dots e_n$ be an orthonormal basis for $H^1_{odd}(\tilde{X},
\reals)$ with respect to $\langle \cdot, \cdot \rangle$, so that
\begin{displaymath}
A^*A e_i = \lambda_i^2 e_i, 
\end{displaymath}
where $e^t = \lambda_1 \ge \dots \ge \lambda_n = e^{-t} > 0$. Then
$\|A e_i \| = \lambda_i$. Note that $e_1$ is the imaginary part of the
holomorphic $1$-form $\tilde{\omega}$ on $\tilde{X}$ for which
$\tilde{\omega}^2 = q$, and $e_n$ is the real part of
$\tilde{\omega}$. We can read off that $\lambda_1 = e^{t}$, $\lambda_n
= e^{-t}$. Also since $A$ is symplectic, $\lambda_2 =
\lambda_{n-1}^{-1}$. 

Note since we are assuming $q \in K$, 
the pair $(\tilde{X},\tilde{\omega})$
associated to $q$ belongs to a compact subset (depending only on $K$)
of the space of Abelian differentials. 
Then, by Theorem~\ref{theorem:forni},
\begin{equation}
\label{eq:lambda1:crude}
\lambda_2  =  O(e^{t}) \text{ and } \lambda_{n-1}^{-1} = O(e^{t}),
\end{equation}
and if (\ref{eq:s:half:time}) holds,
\begin{equation}
\label{eq:lambda1:fine}
\lambda_{2} =  O(e^{(1 - \alpha)t }) \text{ and } \lambda_{n-1}^{-1} =
O(e^{(1-\alpha)t}). 
\end{equation}
Since $A$ has determinant $1$, the product of the
$\lambda_i$ is $1$.

We can identify the tangent space to $\mathcal{Q} \cT_g$ with 
$H^1_{odd}(\tilde{X},\reals^2) = H^1_{odd}(\tilde{X},\reals) \tensor
\reals^2$. 

Let $e_i^+ = e_i \tensor \begin{pmatrix} 1 \\ 0 \end{pmatrix}$ and
$e_i^- = e_i \tensor \begin{pmatrix} 0 \\ 1 \end{pmatrix}$ so that the
$e_i^+$ 
and $e_i^-$ together form a basis for the tangent space
to $\mathcal{Q} \cT_g$ at $q$.
Then $(dg_t)_*(e_i^+) = e^{t} A e_i^+$, and $(dg_t)_*(e_i^-) = e^{-t}
A e_i^-$. (In the above we adopt the convention that
$A$ acts on the tensor product 
$H^1_{odd}(\tilde{X},\reals) \tensor \reals^2$ by acting on the 
$H^1_{odd}(\tilde{X},\reals)$ factor, and trivially on the second
factor.)

Let $u_1 = e_1^+-e_n^-$. We can complete the set
$\{u_1\}$ to a basis $\{ u_1, \dots, u_{n-1} \}$ for the tangent space
of $\qonetg$ at $q$, so that for $2 \le i \le n-1$, $u_i = e_i^+ + w_i$,
where $w_i$ is in the span of $\{e_2^-,\dots, e_{n-1}^-\}$. For
consistency, we let $w_1 =- e_n^-$. 
Note that the $\|w_i\|$ are bounded since $q$ belongs to a 
compact set $K$, and the leaves of $\cF^{u}$ are transverse to the
spheres $S(X)$ on $K$. 

We extend the Hodge inner product $\langle \cdot, \cdot \rangle$
to $H^1_{odd}(\tilde{X},\reals^2)$ by declaring that $\langle
e_i^+,e_j^-\rangle = 0$ for all $i,j$. This also extends the Hodge
norm. We now estimate
\begin{align}
\label{eq:wedge:norm}
\| (dg_t)_* (u_1 \wedge \dots \wedge u_{n-1}) \|^2 & = 
\langle u_1 \wedge\dots\wedge u_{n-1}, [(dg_t)_*]^* [(dg_t)_*]
u_1 \wedge\dots\wedge u_{n-1} \rangle \\
& = \sum_{P} \langle u_1\wedge\dots\wedge u_{n-1}, 
[(dg_t)_*]^* [(dg_t)_*] (\bigwedge_{i \not\in P} e_i^+ \wedge
\bigwedge_{i \in P} w_i) \rangle \notag \\
& = \sum_{P} \langle u_1\wedge\dots\wedge u_{n-1}, 
(\bigwedge_{i \not\in P} e^{2t} \lambda_i^2 e_i^+ \wedge \bigwedge_{i
  \in   P} e^{-2t} A^*A w_i) \rangle \notag 
\end{align}
where $[(dg_t)_*]^*$ is the adjoint of the linear transformation
$(dg_t)_*$, and the sum is over all subsets $P$ of $\{1,\dots,n-1\}$. 

Note that by (\ref{eq:lambda1:crude}), for $1 \le i \le n-1$
\begin{displaymath}
\| e^{-2t} A^*A w_i \| \le C \le C \| e^{2t} \lambda_i^2 \|.
\end{displaymath}
Thus all the terms in (\ref{eq:wedge:norm}) are bounded by
$C e^{2t(n-1)} \lambda_1^2 \prod_{i=2}^{n-1} \lambda_i^2 = C e^{2 nt} =
C e^{2 ht}$. Thus,
\begin{equation}
\label{eq:crude:normwedge}
\| (dg_t)_* (u_1 \wedge \dots \wedge u_{n-1}) \| = O(e^{ht}). 
\end{equation}
By definition, $\Delta(q,t)$ is
the determinant of the linear transformation $D = (d\pi_Y)_* (dg_t)_*$.
Since $Y$ is in a compact set, the norm of $(d\pi_Y)_*$ is bounded.
This together with (\ref{eq:crude:normwedge}) implies
(\ref{eq:measure:polar:crude:bound}). If we assume
(\ref{eq:lambda1:fine}), we get for any $1 \le i \le n-1$, 
\begin{displaymath}
\|e^{-2t} A^*A w_i \| \le C e^{-2\alpha t} \le C e^{-2\alpha t} 
e^{2t} \lambda_i^2,
\end{displaymath}
which implies that the contribution of all the terms in
(\ref{eq:wedge:norm}) to $\Delta(q,t)$ except for the term 
with $P = \emptyset$ is $O(e^{(h-\alpha)t})$.
It remains to evaluate the term with $P = \emptyset$. 
Dropping the rest of the terms is equivalent to replacing the
map $\Phi(q,t)$ by the composition of three maps:
\begin{enumerate}
\item Projecting from the sphere $S(X)$ to the leaf of the
expanding foliation $\cF^{u}$ passing through $q$. 
\item Flowing by time $t$.
\item Applying the map $\pi_Y$ to project from the leaf of $\cF^{u}$
  through $g_t q$ to Teichm\"uller space $\cT_g$. 
\end{enumerate}
The Jacobian of the resulting map is the product of the
Jacobians from steps 1), 2) and 3). Step 1) gives a factor
of $\lambda^-(q)$, step 2) gives a factor of $e^{ht}$, and step 3)
gives a factor of $\lambda^+(g_t q)$, so one gets formula
(\ref{eq:measure:polar:asymptiotics}). 
\qed\medskip

\subsection{The Hodge Norm and Quadratic Differentials.}
\label{sec:hodge:quadratic}

\subsubsection{The Hodge Norm and the Flat Metric.}\label{sec:hodge:flat}
Theorem~\ref{theorem:hodge:hyperbolic} allows us to estimate the 
Hodge norm in terms of the hyperbolic metric. In this subsection 
we state a lemma which allows us to estimate the Hodge norm in terms
of the flat metric. This is done by first estimating the hyperbolic
metric in terms of the flat metric, and then using
Theorem~\ref{theorem:hodge:hyperbolic}. 

Let $q$ be a holomorphic quadratic differential, and let $\ell_q$
denote length in the flat metric defined by $q$. 
For each $\epsilon > 0$ let $K(\epsilon)$ denote the complex
generated by all saddle connections shorter then $\epsilon$
(see \cite[\S{6}]{Eskin:Masur}). 
Recall that $K(\epsilon)$ contains all saddle
connections shorter then $\epsilon$, and that $\ell_q(\partial
K(\epsilon)) = O(\epsilon)$, where the implied constant depends only
on the genus. (Also if $\mathcal{C}$ is a cylinder whose boundary is
shorter then 
$\epsilon$, we include $\mathcal{C}$ in $K(\epsilon)$). 
Let $\epsilon_1 <  \dots < \epsilon_n$ denote the values of
$\epsilon$ where $K(\epsilon)$ changes. (Note that $n$ is bounded in
terms of the genus). Pick a constant $C \GG 1$ and drop all $\epsilon_i$
such that $\epsilon_i > \epsilon_{i+1}/C$. After renumbering, we
obtain new $\epsilon_i$ $1 \le i \le m$, 
with $\epsilon_i < \epsilon_{i+1}/C$. To simplify notation, we let
$\epsilon_{m+1} = 1$. For $1 \le i \le m$, let
$K_i$ denote the complex $K(\epsilon_i)$. 

\begin{lemma}[Rafi, \cite{Rafi:TT}]
\label{lemma:flat:to:hyperbolic}
Let $q$ be a holomorphic quadratic differential on a surface $\Sigma_{g}$, and
let $\sigma$ be the hyperbolic metric on $\Sigma_{g}$ in the same conformal
class as $q$. Choose $C \GG 1$, and
let the $\epsilon_i$, $1 \le i \le m$ and the complexes $K_i$ 
be as in the previous paragraph. Then
\begin{itemize}
\item[{\rm (a)}] There exists a constant $\epsilon_0$ depending only on $C$ and the genus
(with $\epsilon_0 \to 0$ as $C \to \infty$) such that any simple 
closed curve $\alpha$ on $\Sigma_{g}$ with the hyperbolic length
$\ell_\alpha(\sigma) < \epsilon_0$ 
is homotopic to a connected component of 
the boundary of one of the $K_i$. 

\item[{\rm (b)}] Let $\gamma$ be a connected component of $\partial
  K_i$. If $\gamma$ is not a core curve of a flat annulus, then
$\ell_\gamma(\sigma) \ge C_1 /|\log \epsilon_{i}|$,
where $C_1$ is a constant depending only on the genus. 
\end{itemize}
\end{lemma}

\subsubsection{The Modified Hodge Norm and Quadratic Differentials.}
\label{sec:modified:norm}
Theorem~\ref{theorem:forni} gives a partial hyperbolicity property of
the geodesic flow on spaces of Abelian differentials. In our
applications, we need a similar property for the spaces $\qonetg$ of
quadratic differentials. A standard construction, given $X \in
\cT_g$ and $q$ a holomorphic quadratic differential on $X$, is to pass
to the possibly ramified double cover on which the foliation defined by $q$ is
orientable. This yields a surface $\tilde{X}$ with an holomorphic
Abelian differential $\omega$. However, a major difficulty is the
following: even if $q$ belongs to a compact subset of $\qonetg$, the
complex structure on $\tilde{X}$ may have very short closed curves in
the hyperbolic metric. This may occur since even if one restricts $q$
to compact subsets of Teichm\"uller space, the flat structure defined by
$q$ may have arbitrarily short saddle connections (connecting distinct
zeroes). Such a saddle connection may lift to a very short loop when
one takes a double cover. Other more complicated types of degeneration
are also possible. However, the following simple observation is key to
our approach:

\begin{lemma}
\label{lemma:no:flat:cylinders}
For every compact subset $\cK \subset \cM_g$ there exists a constant
$\epsilon_0 > 0$ depending only on $\cK$ such that for any $q \in
\pi^{-1}(\cK)$, the flat structure associated to the holomorphic
quadratic differential  on the orienting cover of $q$ has no closed
trajectories of Euclidean length less than $\epsilon_0$ which are part
of a flat cylinder. 
\end{lemma}

Together with Lemma~\ref{lemma:flat:to:hyperbolic} (b),
Lemma~\ref{lemma:no:flat:cylinders} will allow us to estimate the
hyperbolic length of short curves on $\tilde{X}$.

\bold{Proof of Lemma~\ref{lemma:no:flat:cylinders}.} Given the set $\cK$, there exists a constant $\epsilon_0 > 0$ such  that any loop of length less than $\epsilon_0$ in the flat metric
is contractible. Now suppose $q \in \pi^{-1}(\cK)$, and 
$\gamma$ is a trajectory on the associated orienting double 
cover which has length less than $\epsilon_0$, and is part of a flat
cylinder $C$. Let $\gamma_1$ be the projection of $\gamma$ 
to the original flat structure defined by $q$. Then $\gamma_1$ must also be a
part of a flat cylinder, and the length of $\gamma_1$ must be at most
$\epsilon_0$. Then, from the definition of $\epsilon_0$ it follows
that $\gamma_1$ must be contractible. But this is impossible since the
curvature of $\gamma_1$ is zero, and we are not allowing $q$ to have
poles.
\qed\medskip

\bold{Short bases.} Suppose $(X,\omega) \in \omegaonetg$ (where the
notation $\omegaonetg$ is defined in \S\ref{sec:hodge:geod}). 
Fix $\epsilon_1 <
\epsilon_*$ (where $\epsilon_*$ is the Margulis constant defined in
\S\ref{sec:hodge:hyp})
and let $\alpha_1, \dots, \alpha_k$ be the curves with
hyperbolic length less than $\epsilon_1$ on $X$. For $1 \le i \le k$,
let $\beta_i$ be the shortest curve in the flat metric defined by
$\omega$ with $i(\alpha_i, \beta_i) =1$. Let $\gamma_r$, $1 \le r \le
2g-2k$ be moderate length curves on $X$ so that the $\alpha_j$,
$\beta_j$ and $\gamma_j$ are a symplectic basis $\cS$ for
$H_1(X,\reals)$. We will call such a basis {\em short.}

\bold{The functions $\chi_i(\omega)$ and the modified Hodge norm.}
We would like to use the Hodge norm on double covers of surfaces 
in $\qonetg$. 
One difficulty is that if one takes the double cover
of some quadratic differential in the multiple zero locus, the Hodge
norm in some directions tangent to the multiple zero locus might
vanish. As a consequence, if we use the
Hodge norm to define a ``Hodge distance'' on $\qonetg$, then we find
that some balls are not compact, and the resulting ``distance'' does
not separate points on the multiple zero locus. 

We now define a modification of the Hodge norm in order to avoid these
problems. The modified norm is defined on the tangent space to the
space of pairs $(X,\omega)$ where $X$ is a Riemann surface and
$\omega$ is a holomorphic $1$-form on $X$. Unlike the Hodge norm, 
the modified Hodge norm will
depend not only on the complex structure on $X$ but also on the choice
of a holomorphic $1$-form $\omega$ on $X$. Let $\{\alpha_i, \beta_i,
\gamma_r\}_{1 \le i \le k, 1 \le r \le 2g-2k}$ 
be a short basis for $(X,\omega)$. 
For $1 \le i \le k$ we define $\chi_i(\omega)$
to be $0$ if $\alpha_i$ has a flat annulus in the flat metric defined
by $\omega$ and $1$ otherwise. 

We can write any $\theta \in H^1(X,\reals)$ as
\begin{equation}
\label{eq:expand:in:basis}
\theta = \sum_{i=1}^k a_i (*c_{\alpha_i}) + \sum_{i=1}^k b_i
\ell_{\alpha_i}(\sigma)^{1/2} (*c_{\beta_i})  + \sum_{r=1}^{2g -2k } u_r (*c_{\gamma_r}).
\end{equation}
(Recall that $\sigma$ denotes the hyperbolic metric in the conformal
class of $X$, and for a curve $\alpha$ on $X$, $\ell_\alpha(\sigma)$
denotes the length of $\alpha$ in the metric $\sigma$). 
We then define
\begin{equation}
\label{eq:def:modified:hodge:norm}
 \|\theta\|' = \|\theta\| +  \left( \sum_{i=1}^k \chi_i(\omega) |a_i| +
   \sum_{i=1}^k |b_i| +
    \sum_{r=1}^{2g -2k} |u_r| \right).
\end{equation}
From (\ref{eq:def:modified:hodge:norm}) we have for $1 \le i \le k$, 
\begin{equation}
\label{eq:star:c:alpha:prime:norm}
\|*\!c_{\alpha_i}\|' \approx 1, 
\end{equation}
as long as $\alpha_i$ has no flat annulus in the metric defined by
$\omega$. Similarly, from (\ref{eq:c:beta:holomorphic})
we have 
\begin{equation}
\label{eq:star:c:beta:prime:norm}
\|*\!c_{\beta_i}\|' \approx  \| *\!c_{\beta_i} \| \approx \frac{1}{\ell_{\alpha_i}(\sigma)^{1/2}}.
\end{equation}
In addition, in view of Theorem~\ref{theorem:hodge:hyperbolic}, 
if $\gamma$ is any other moderate length curve on $X$, 
$\|*\!c_\gamma\|' \approx \|*\!c_\gamma \| = O(1)$. 
Thus, if $\cB$ is a short basis associated to $\omega$, then 
for any $\gamma \in \cB$, 
\begin{equation}
\label{eq:short:basis:extremal:length}
\Ext_\gamma(\omega)^{1/2} \approx \|\!*\!c_\gamma\| \le
\|\!*\!c_\gamma\|' \end{equation}
(cf. (\ref{eq:extremal:hodge})). (By $\Ext_\gamma(\omega)$ we mean the
extremal length of $\gamma$ in $X$, the conformal structure defined 
by $\omega$). 

\bold{Remark.} From the construction, we see that the modified Hodge
norm is greater than the Hodge norm. Also, if the flat length of 
shortest curve in
the flat metric defined by $\omega$ is greater than $\epsilon_1$, then
for any cohomology class $\lambda$, for some $C$ depending on
$\epsilon_1$ and the 
genus, 
\begin{equation}
\label{eq:modified:hodge:compare:to:hodge}
\|\lambda\|' \le C \|\lambda\|,
\end{equation}
i.e. the modified Hodge norm is within a multiplicative constant of
the Hodge norm.

\subsection{Non-expansion of the modified Hodge norm.}
In this subsection we will prove a weak version of
Theorem~\ref{theorem:forni} for the modified Hodge norm.

Let $\| \cdot \|'_0$ denote the modified Hodge norm on the surface
$S_0 = S = (X_0, \omega_0)$, and let 
$\| \cdot \|'_t$ denote the modified Hodge norm on the surface $S_t =
g_t S = (X_t,\omega_t)$.
\begin{theorem}
\label{theorem:nonincrease:modified:norm}
There exists a constant $C$ depending only on the genus and on
$\epsilon_1$, such that 
for any $\theta \in H^1(X,\reals)$ and any $t \ge 0$, 
\begin{equation}
\label{eq:nonincrease:modified:norm}
\|\theta\|'_t \le C e^{t} \| \theta\|'_0.
\end{equation}
\end{theorem}
\medskip

The proof of Theorem~\ref{theorem:nonincrease:modified:norm} 
is based on the following lemma:
\begin{lemma}[Rafi]
\label{lemma:rafi:short}
Suppose $(X_0,\omega_0) \in \omegaonetg$, and write $(X_t,\omega_t)$ for
$g_t (X_0,\omega_0)$. Suppose $\cB$ is a short basis in the flat
structure defined by $\omega_0$, and $\cB'$ is a short basis in the
flat structure defined by $\omega_t$. Suppose $\alpha$ is a curve of
hyperbolic length $< \epsilon_1$ \mc{check constants} on both $X_0$
and $X_t$ (so that in particular $\alpha \in \cB$ and $\alpha \in
\cB'$), and that $\alpha$ has no flat annulus. Let $\beta' \in \cB'$
denote the curve with $i(\alpha,\beta') = 1$. Then, for any $\gamma
\in \cB$, we have
\begin{displaymath}
i(\gamma,\beta') \le C e^t \Ext_\gamma(X_0)^{1/2},
\end{displaymath}
where $i(\cdot,\cdot)$ denotes the intersection number.
\end{lemma}

Recall that given simple closed curves $\alpha$ and $\beta$ on $\Sigma_g$, the
intersection number $i(\alpha,\beta)$ is the minimum number of points in
which representatives of $\alpha$ and $\beta$ must intersect.

\bold{Proof of Lemma~\ref{lemma:rafi:short}.}  First we claim that $i(\gamma, \beta')\leq i(\gamma,
P)$, where $P$ is the shortest pants decomposition of $X_t$. 
This is
because, either $\gamma$ is disjoint from $\alpha$ or if it intersects
$\alpha$, the relative twisting of $\gamma$ and $\beta'$ is bounded
(see \cite[\S{4}]{Rafi:CM}). But this shortest pants
decomposition has extremal length less than a constant times $e^{2t}$
in $X_0$.  Using \cite[Lemma 5.1]{Minsky:geodesics} we have
$$ 
i(\gamma, P)^2 \le \Ext_\gamma(X_0) \Ext_P(X_0) 
\leq C^2 e^{2t} \Ext_\gamma(X_0).
$$ 
This finishes the proof. 
\qed\medskip

\bold{Proof of Theorem~\ref{theorem:nonincrease:modified:norm}.}
Let $\alpha_1', \dots, \alpha_k'$ be the curves with hyperbolic
length less than $\epsilon_1$ on $X_t = \pi(S_t)$.
Let $\beta_1', \dots \beta_k'$
be the shortest curves in the Euclidean metric on $S_t$ 
such that $i(\alpha_r, \beta_j) = \delta_{rj}$. 
Let $\gamma_r'$, $1 \le r \le 2g-k$ be moderate length curves on $X_t$
so that the $\alpha_j$, $\beta_j'$ and $\gamma_j'$ are a symplectic
basis $\cS'$ for $H_1(X,\reals)$. Thus $\cS'$ is a short basis. 
For any $\theta \in
H^1(X,\reals)$, we may write
\begin{equation}
\label{eq:expand:in:prime:basis}
\theta = \sum_{i=1}^k a_i (*c_{\alpha_i'}) + \sum_{i=1}^k b_i
\ell_{\alpha_i'}(\sigma)^{1/2} (*c_{\beta_i'})  + \sum_{r=1}^{2g -k } u_i (*c_{\gamma_r'}).
\end{equation}
In view of (\ref{eq:def:modified:hodge:norm}) and
Theorem~\ref{theorem:forni}, it is enough to show that
\begin{equation}
\label{eq:bound:coeffs}
\sum_{i=1}^k \chi_i(\omega_t) |a_i| +
   \sum_{i=1}^k |b_i| +
    \sum_{r=1}^{2g -k} |u_r|  \le C e^{t} \|\theta\|'_0.
\end{equation}
We have
\begin{displaymath}
b_i = \frac{1}{\ell_{\alpha_i'}(\sigma)^{1/2}} \langle c_{\alpha_i'}, \theta
  \rangle_t \le   \frac{1}{\ell_{\alpha_i'}(\sigma)^{1/2}} \|
  c_{\alpha_i'}\|_t \|\theta \|_t \le \| \theta\|_t \le C e^{t}
  \|\theta\|_0 \le C e^{t} \|\theta\|'_0.
\end{displaymath}
Here, $\langle \cdot, \cdot \rangle_t$ denotes the Hodge inner product
on $X_t$.  
Also let $\delta_i' \in \cS'$ be the element with
$i(\gamma_i',\delta_i') = 1$. Then,
\begin{displaymath}
|u_i| = \langle c_{\delta_i'}, \theta
  \rangle_t \le   \| 
  c_{\delta_i'}\|_t \|\theta \|_t \le C_1 \| \theta\|_t \le C_2 e^{t}
  \|\theta\|_0 \le C_3 e^{t} \|\theta\|'_0,
\end{displaymath}
where $C_1$, $C_2$ and $C_3$ depend only on the genus. It remains to bound
$|a_i|$. We may assume that  
$\alpha_i$ has length less than $\epsilon_1$ on $X_0$
as well. (If $\alpha_i$ is longer then $\epsilon_1$ on $X_0$, we can
put in an intermediate point $X_0'$ where $\alpha_i$ has length
$\approx \epsilon_1$.) \mc{explain more}

Let $\cS = \{\alpha_i,\beta_i,\gamma_i \}$ be a short basis for
$S_0$. It is enough to prove (\ref{eq:bound:coeffs}) assuming $\theta
= *c_\lambda$ for some $\lambda \in \cS$. 
Then, 
\begin{displaymath}
a_i = \langle *c_{\beta_i'}, c_\lambda \rangle = I(\lambda,\beta_i').
\end{displaymath}
Since the algebraic intersection number $I(\cdot,\cdot)$ is bounded by
the geometric intersection number $i(\cdot,\cdot)$, we have
\begin{displaymath}
|a_i| = |I(\lambda,\beta_i')| \le i(\lambda,\beta_i') \le C e^t
\Ext_\lambda(X_0)^{1/2} \le C e^t \|\!*\!c_\lambda\|',
\end{displaymath}
where we have used Lemma~\ref{lemma:rafi:short} and
(\ref{eq:short:basis:extremal:length}). 
\qed\medskip

\subsection{The Euclidean, Teichm\"uller and Hodge Distances.}
\label{sec:euclidean:hodge}

Note that we can locally identify a leaf of $\cF^s$ (or $\cF^u$) with
a subspace of $H^1(\tilde{X},\reals)$, where $\tilde{X}$ is the double
cover of $X \in \cT_g$.  
If $\gamma$ is a map from $[0,r]$ into some leaf of $\cF^{ss}$, then
we define the (modified) Hodge length $\ell(\gamma)$ of $\gamma$ as $\int_0^r
\|\gamma'(t)\|' \, dt$, where $\| \cdot \|'$ is the modified Hodge
norm. If $q$ and $q'$ belong to the same leaf of $\cF^{ss}$, then we define
the Hodge distance
$d_H(q,q')$ to be the infimum of $\ell(\gamma)$ where $\gamma$ varies over
paths connecting $q$ and $q'$ and staying in the leaf of $\cF^{ss}$
containing $q$ and $q'$. We make the same definition if $q$ and $q'$
are on the same leaf of $\cF^{uu}$.

In what follows $\cK$ is a compact subset of $\cT_g$, 
and all implied constants depend on
$\cK$. Given the set $\cK$, there exists a constant $\epsilon_0 > 0$ such 
that for any $q \in \pi^{-1}(\cK)$ 
any loop of length less than $\epsilon_0$ in the flat metric 
given by $q$ is contractible.

For each $q \in \pi^{-1}(\cK)$ there exists a canonical (marked) Delaunay
triangulation of the flat metric given by $q$ (see
\cite[\S{4}]{Masur:Smillie}). Since there are 
only finitely many combinatorial types of (unmarked) triangulations
on a genus $g$ surface, and the mapping class group $\Gamma$ acts
properly discontinuously, there are only finitely many combinatorial
types $\alpha$ of marked Delaunay triangulations on $\pi^{-1}(\cK)$. Call this
set $I$. 
We have
\begin{displaymath}
\pi^{-1}(\cK) =  \bigsqcup_{\alpha \in I} W_\alpha,
\end{displaymath}
where $W_\alpha$ is the subset of $\pi^{-1}(\cK)$ where the combinatorial type of
the marked Delaunay triangulation is given by $\alpha$.
Let $J \subset I$ denote the subset such that for $\alpha \in J$,
$W_\alpha$ is relatively open in $\pi^{-1}(\cK)$. Then we may write
\begin{displaymath}
\pi^{-1}(\cK) =  \bigsqcup_{\alpha \in J} W_\alpha',
\end{displaymath}
where $W_\alpha \subseteq W_\alpha' \subseteq \overline{W_\alpha}$.
If we choose a basis $\cS$ for $H_1^{odd}(\tilde{X},\zed)$, then we
have a period map $\Phi_{\cS}: \qonetg \to \cx^h$ given by integrating
the canonical square root of $q$ along the chosen basis. Recall that
this is a local coordinate system for $\qonetg$. In fact the following
holds:

\begin{lemma}
\label{lemma:good:basis}
For each $\alpha \in J$ there exists a basis $\cS_\alpha$ of
$H_1^{odd}(\tilde{X},\zed)$ consisting of lifts of edges of the
Delaunay triangulation such that the map $\Phi_\alpha =
\Phi_{\cS_\alpha}$ is defined on all of $W_\alpha'$
  and is a coordinate system on the part
  of $W_\alpha'$ not contained in the multiple zero locus. 
\end{lemma}

\bold{Proof.} See \cite[\S{4}]{Masur:Smillie}. 
\qed\medskip

We can now define a ``Euclidean norm'' on the tangent space of $\qonetg$
as follows: if $v$ is a tangent vector based at the point $q \in
W_\alpha$ we define
\begin{displaymath}
\|v\|_E = | \Phi_{\alpha}^*(v) |,
\end{displaymath}
where $| \cdot |$ is the standard norm on $\cx^h$. 
The norm $\| \cdot \|_E$ depends on the choices of the $\cS_\alpha$
and the $W_\alpha'$;
however if $\| \cdot \|'_E$ is the Euclidean norm obtained from a
different choices, then the ratio of $\| \cdot \|_E$ and
$\|\cdot \|_E'$ is bounded by a constant depending only on $\cK$. 

\bold{The Euclidean Distance.}
We define the Euclidean length of a path $\gamma: [0,r] \to \qonetg$
to be $\int_0^r \| \gamma'(t)\|_E \, dt$. We define the Euclidean distance 
$d_E(q,q')$ to be the infimum of the Euclidean length of paths
connecting $q$ and $q'$. Locally, up to a multiplicative constant, 
$d_E(q, q') = | \Phi_\alpha(q) - \Phi_\alpha(q')|$, 

We denote the Teichm\"uller distance on $\cT_g$ by
$d_\cT$. 
For $p \in \qonetg$ and $\epsilon
>0$ let $B_E(p,\epsilon)$ denote the ball in the Euclidean metric
of radius $\epsilon$ centered at $p$, and let 
$B_{\cT}(\pi(p), \epsilon)$ denote the analogous ball in the Teichm\"uller metric. We note the following:
\begin{lemma}
\label{lemma:Teich:equivalent:to:euc}
There exist continuous functions $f_i: \reals^+ \to \reals^+$ with $f_i(r)
\to 0$ as $r \to 0$ such that for any two points $p_1$ and $p_2$ in
$\pi^{-1}(\cK)$ on the same leaf of $\cF^{uu}$, we have
\begin{displaymath}
f_1(d_E(p_1,p_2)) \le d_{\cT}(\pi(p_1),\pi(p_2)) \le f_2(d_E(p_1,p_2)).
\end{displaymath}
\end{lemma}

\bold{Proof.} This follows from the compactness of $\cK$. 
\qed\medskip

\bold{Remark.} The Euclidean distance has a number of useful
properties: it behaves well near the multiple zero locus, and on
compact subsets of $\cT_g$ it is
absolutely continuous with respect to the Teichm\"uller distance (see
Lemma~\ref{lemma:Teich:equivalent:to:euc}). However, often we need to deal
instead with the modified Hodge distance because of its non-expansion and
decay properties under the geodesic flow
(Theorem~\ref{theorem:hodge:distance:nonincrease} below). 

The following theorem is the main result of this subsection:
\begin{theorem}
\label{theorem:hodge:equivalent:to:euclidean}
There exists constants $0 < c_1 < c_2$, both depending only on $\cK$, 
such that for any $q_1, q_2 \in \pi^{-1}(\cK)$ on the same leaf of
$\cF^{uu}$ with
$d_E(q_1,q_2) < 1$ we have:
\begin{displaymath}
c_1 d_E(q_1, q_2) \le d_H(q_1,q_2) \le c_2 d_E(q_1,q_2) |\log d_E(q_1,q_2)|^{1/2}.
\end{displaymath}
\end{theorem}
Thus, in particular, on compact sets, 
the modified Hodge metric is equivalent to the Euclidean metric, and
then, in view of Lemma~\ref{lemma:Teich:equivalent:to:euc}, the
modified Hodge metric on the restriction of a leaf of $\cF^{uu}$ to
$\pi^{-1}(\cK)$ is equivalent to the Teichm\"uller metric.

Note that the lower bound in
Theorem~\ref{theorem:hodge:equivalent:to:euclidean} is clear since up
to a constant depending only on $\cK$, 
the modified Hodge norm is always
bigger than the Euclidean norm, see (\ref{eq:star:c:alpha:prime:norm})
and (\ref{eq:star:c:beta:prime:norm}).
The rest of this subsection will consist of the proof of the upper
bound in Theorem~\ref{theorem:hodge:equivalent:to:euclidean}. 
For $q \in \qonetg$, let $\ell_{min}(q)$ denote the length of the shortest
saddle connection in the flat metric defined by $q$.

We will use the following standard fact about Delaunay triangulations:
\begin{lemma}
\label{lemma:delaunay}
Suppose $T$ is the Delaunay triangulation of a surface in $q \in
\qonetg$. Let $\ell_{min}(q)$ denote the length of the shortest saddle
connection in $q$. Then any saddle connection in $q$ of length at most
$\sqrt{2} \ell_{min}(q)$ is an edge of $T$.
\end{lemma}

\bold{Proof.} Let $z_1, \dots, z_n$ denote the zeroes of $q$ (i.e. the
conical points in the flat metric). For each $z_i$, let the {\em
Voronoi polygon} $V_i$ denote the set of points of $q$ that are closer to
$z_i$ then to any other $z_j$, $j \ne i$. 

Recall that the Delaunay triangulation $T$ of $q$ is dual to the Voronoi
diagram in the sense that $z_i$ and $z_j$ are connected by an edge of 
$T$ if $V_i$ and $V_j$ share a common edge (and also under some
conditions if $V_i$ and $V_j$ share a common vertex). 

Let $e$ be a saddle connection in $q$,
connecting $z_i$ and $z_j$. (See Figure~\ref{fig:delaunay}). Let $m$
be the midpoint of $e$. Let $k$ be such that $d(z_k,m)$ is minimal. 
Suppose $e$ is {\em not} an edge of the Delaunay triangulation
$T$. Then $k \ne i,j$ and $d(m,z_k) \le d(m,z_i) = \frac{1}{2}
d(z_i,z_j)$. We may assume that the angle at $m$ between the segments
$\overline{m z_i}$ and $\overline{m z_k}$ is smaller then $\pi/2$
(otherwise replace $z_i$ by $z_j$). Let $p$ be the point on the
segment $\overline{m z_i}$ such that $d(p,m) = d(p, z_k)$. 
Now consider the isosceles acute triangle $\tau$
whose vertices are $m$, $p$ and $z_k$. This triangle cannot contain
any zeroes (or else $z_k$ would not be the saddle connection closest
to $m$). Therefore, $d(p,z_k) \le \sqrt{2}d(m,z_k)$. Hence,
\begin{multline*}
\ell_{min}(q) \le d(z_i, z_k) \le d(z_i,p) + d(p,z_k) \le (d(z_i,m) -
d(m,z_k))+\sqrt{2}d(m,z_k) \\
= \frac{1}{2} d(z_i, z_j) + (\sqrt{2}-1) d(m,z_k) \le 
\frac{1}{2} d(z_i, z_j) + \frac{(\sqrt{2}-1)}{2} d(z_i,z_j) = \frac{\sqrt{2}}{2}d(z_i,z_j).
\end{multline*}
Thus, if $e$ is not an edge of the Delaunay triangulation, then 
$\ell(e) \ge \sqrt{2} \ell_{min}(q)$. 
\qed\medskip

\makefig{Proof of Lemma~\ref{lemma:delaunay}.}{fig:delaunay}{\begin{picture}(0,0)%
\includegraphics{delaunay.pstex}%
\end{picture}%
\setlength{\unitlength}{4144sp}%
\begingroup\makeatletter\ifx\SetFigFont\undefined%
\gdef\SetFigFont#1#2#3#4#5{%
  \reset@font\fontsize{#1}{#2pt}%
  \fontfamily{#3}\fontseries{#4}\fontshape{#5}%
  \selectfont}%
\fi\endgroup%
\begin{picture}(1875,2220)(4981,-5251)
\put(4996,-5191){\makebox(0,0)[lb]{\smash{{\SetFigFont{10}{12.0}{\rmdefault}{\mddefault}{\updefault}{\color[rgb]{0,0,0}$z_i$}%
}}}}
\put(6841,-3166){\makebox(0,0)[lb]{\smash{{\SetFigFont{10}{12.0}{\rmdefault}{\mddefault}{\updefault}{\color[rgb]{0,0,0}$z_j$}%
}}}}
\put(6301,-4561){\makebox(0,0)[lb]{\smash{{\SetFigFont{10}{12.0}{\rmdefault}{\mddefault}{\updefault}{\color[rgb]{0,0,0}$z_k$}%
}}}}
\put(5716,-4066){\makebox(0,0)[lb]{\smash{{\SetFigFont{10}{12.0}{\rmdefault}{\mddefault}{\updefault}{\color[rgb]{0,0,0}$m$}%
}}}}
\put(5446,-4471){\makebox(0,0)[lb]{\smash{{\SetFigFont{10}{12.0}{\rmdefault}{\mddefault}{\updefault}{\color[rgb]{0,0,0}$p$}%
}}}}
\end{picture}%
}

Recall that an integral multicurve is a finite set of
oriented simple closed curves with integral weights. By convention, a negative
weight corresponds to reversing the orientation. 

\begin{lemma}
\label{lemma:goodmulticurve}
Suppose $T$ is a geodesic triangulation of an orientable
surface $\omega$ in $\omegaonetg$. 
We orient each edge of $T$ so that its $x$ component is
positive (for a vertical edge the orientation is not defined). 
Suppose $W$ is a subset of the edges of $T$. Then there
exists an integral multicurve $\Delta$ on $\omega$ such that
\begin{itemize}
\item[{\rm (a)}] $\Delta$ is disjoint from the vertices of $T$ and is
  transverse to the edges.
\item[{\rm (b)}] $\Delta$ crosses each edge of $W$ at least once.
\item[{\rm (c)}] Each time  $\Delta$ crosses an  
edge of $e \in W$, 
the crossing is from left to right (with respect to the
orientation on the edge). If $e$ is vertical, then all crossings
must be from the same side. 
\item[{\rm (d)}] For each edge $e$ of $T$, the intersection number
$i(\Delta,e) \le n$, where $n$ depends only on the genus $g$. 
\end{itemize}
\end{lemma}

\makefignocenter{Lemma~\ref{lemma:goodmulticurve}.
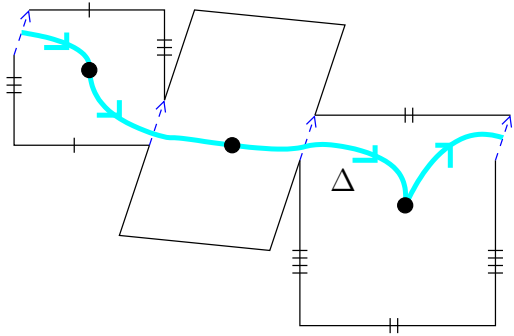
The surface consists of three flat tori glued to each other as shown.
The set $W$ consists of the three dotted lines separating the tori.
(The two dotted lines on the opposite sides of the figure are
identified). The multicurve (in this case closed curve) $\Delta$ is 
drawn. The black dots are vertices of the graph $\cG$ (used in the proof).}{fig:graph}{\begin{picture}(0,0)%
\includegraphics{graph.pstex}%
\end{picture}%
\setlength{\unitlength}{4144sp}%
\begingroup\makeatletter\ifx\SetFigFont\undefined%
\gdef\SetFigFont#1#2#3#4#5{%
  \reset@font\fontsize{#1}{#2pt}%
  \fontfamily{#3}\fontseries{#4}\fontshape{#5}%
  \selectfont}%
\fi\endgroup%
\begin{picture}(3045,2004)(5254,-5158)
\put(7201,-4291){\makebox(0,0)[lb]{\smash{{\SetFigFont{12}{14.4}{\rmdefault}{\mddefault}{\updefault}{\color[rgb]{0,0,0}$\Delta$}%
}}}}
\end{picture}%
}

\bold{Proof.}  Each triangulated 
surface which has vertical edges is
the limit of triangulated surfaces which do not. Hence, 
without loss of generality, we may assume that
$T$ has no vertical edges.
 
We define a directed graph $\cG$ as follows (See Figure~\ref{fig:graph}).
The vertices of $\cG$ are the connected components of $\omega \setminus W$. Two
vertices $A$ and $B$ of $\cG$ are connected by an directed edge of
$\cG$ if and only if there exists a saddle connection $\gamma \in W$
such that $A$ is incident to $\gamma$ from the left and $B$ is
incident to $\gamma$ from the right.

We now claim that for each pair $A$, $B$ of vertices of 
$\cG$ there exists at least one directed path from $A$ to
$B$. This can be derived from the minimality of the
  flow in an almost vertical direction, but we prefer to give a direct
  combinatorial argument. Suppose this is not true. Let $\cA$ be the 
set of vertices of $\cG$ 
which can be reached from $A$ by a directed path.
Then $A \in \cA$ and $B \in \cA^c$; in particular $\cA$ and $\cA^c$ are
both non-empty. Now let $D$ be the closure of the connected components
in $\omega \setminus W$ corresponding to vertices of $\cA$. Then $D$ is a 
subsurface of $\omega$, and $\partial D$ consists of edges from $W$. 
For each edge $\gamma \in \partial D$ let
$\epsilon_\gamma = +1$ if $D$ is on the left of $\gamma$ and $-1$
otherwise. Then, (since $D$ is orientable) $\sum_{\gamma \in \partial
  D} \epsilon_\gamma \gamma = 0$. It follows that there exists
$\gamma \in \partial D$ such that the $x$ component of
$\epsilon_\gamma \gamma$ is positive. Since the $x$ component of
$\gamma$ is positive, this implies that $\epsilon_\gamma = 1$. Thus
the directed edge of $\cG$ corresponding to $\gamma$ is directed
from $D$ to $D^c$. This means that some vertex of $\cG$ outside
of $\cA$ can be reached from a vertex in $\cA$ 
by a directed path, contradicting the definition of $\cA$. 

Thus, every pair of vertices of $\cG$ can be connected by a directed
path. This implies every directed edge
of $\cG$ is contained in a directed cycle of $\cG$. Therefore there
exists a finite union $\Delta$ of directed cycles such that every edge
in $\cG$ is contained in $\Delta$. The length of $\Delta$ is bounded
in terms of the size of $\cG$, i.e. in terms of the genus. We can thus
realize $\Delta$ as a multicurve on $\omega$ with properties
(a)-(d). 
\qed\medskip

\begin{lemma}
\label{lemma:expand:one:step}
There exist  $\rho_1 > 0$ and $C > 0$ depending only on $\cK$
such that for any $q \in \pi^{-1}(\cK)$ there
exists $q' \in \pi^{-1}(\cK) \cap \cF^{uu}(q)$ and a path $\gamma: [0,\rho_1]
\to \cF^{uu}(q)$ such that $\gamma(0) = q$, $\gamma(\rho_1) = q'$, 
and for $t \in [0,\rho_1]$, 
\begin{equation}
\label{eq:shortest:increases}
\ell_{min}(\gamma(t)) \ge \ell_{min}(\gamma(0)) \sqrt{1+t^2},
\end{equation}
and also
\begin{equation}
\label{eq:hodge:distance:bound}
d_H(q,q') \le C \int_{0}^{\rho_1} {|\log( \ell_{min}(q) \sqrt{1+t^2})|^{1/2}}{dt}.
\end{equation}
\end{lemma}

\bold{Proof.} Let $\tilde{S} = (\tilde{X},\tilde{\omega})$ denote the
double cover corresponding to $q$. Consider the Delaunay
triangulation of $\tilde{S}$. By
Lemma~\ref{lemma:delaunay}, all saddle connections of length
at most $\sqrt{2} \ell_{min}(\tilde{S})$ belong to the Delaunay
triangulation. Let $W$ denote this set of saddle connections. 
Let $\Delta$ be the multicurve obtained by applying
Lemma~\ref{lemma:goodmulticurve} to $\tilde{S}$ and $W$. 

Let $\tau$ be the involution corresponding to $q$. Let
$\Delta' = \Delta - \tau(\Delta)$. Then $\Delta'$ also
has properties (a)-(d) of Lemma~\ref{lemma:goodmulticurve}. 
In addition, $\tau(\Delta') = -\Delta'$. 

Let $\rho_1 \in (0,\sqrt{2})$ be a constant which will be chosen later
(depending only on the genus). 
We now use $\Delta'$ to define a path $\gamma: [0,\rho_1] \to \qonetg$
contained in $\pi^{-1}(\cK)$ and staying on the same leaf of
$\cF^{uu}$. We will actually
define the path $\tilde{\gamma}(t)$ where for each $t$,
$\tilde{\gamma}(t)$ is the double cover of $\gamma(t)$.
Let $\tilde{S}$ denote the double cover of $q$, and set
$\tilde{\gamma}(0) = \tilde{S}$. 
For each $t \in [0,\rho_1]$, $\tilde{\gamma}(t)$ is
built from the 
same triangles as $\tilde{S}$, but for each edge $e$, 
we add $I(e,\Delta') t \ell_{min}(q)$ 
to the $x$-component, where $I(\cdot, \cdot)$ 
denotes the algebraic intersection number
(cf. \cite[\S{6}]{Masur:Smillie}). 

Note that we do not assert that the Delaunay triangulation of
$\tilde{\gamma}(t)$ is the same as that of $\tilde{S} =
\tilde{\gamma}(0)$. However, because of the form of $\Delta'$, the
involution $\tau$ acts on $\tilde{\gamma}(t)$, and we can let
$\gamma(t)$ be the quotient by $\tau$.

We now claim that (\ref{eq:shortest:increases}) holds. 
Indeed, if $e \in W$ is a saddle connection in
$\tilde{\gamma}(0)$ with 
vector $(x_0,y_0)$, then $x_0 \ge 0$, and on $\tilde{\gamma}(t)$, $e$
has length
\begin{equation}
\label{eq:length:UcapW}
\sqrt{y_0^2 + (x_0 + I(e,\Delta')\ell_{min}(q) t)^2} \ge \sqrt{\ell_{min}(q)^2 +
  \ell_{min}(q)^2   t^2} = \ell_{min}(q) \sqrt{1+t^2}. 
\end{equation}
Suppose $t_1,t_2 \in [0,\rho_1]$ and $\eta$ are such that $\eta$ is a
saddle connection on $\tilde{\gamma}(t)$ for $t_1 < t <t_2$.
Let $\ell_t(\eta)$ denote the length of $\eta$ in the flat metric on
$\tilde{\gamma}(t)$. Then, 
\begin{equation}
\label{eq:change:in:length}
\ell_t(\eta) \ge \ell_{t_1}(\eta) - |I(\eta,\Delta)| \ell_{min}(q) (t-t_1).
\end{equation}

Now suppose $z$ and $w$ be any two zeroes of $q$. 
Let $\lambda_t(z,w)$ denote the shortest path between $z$ and $w$ on
$\tilde{\gamma}(t)$, and let $|\lambda_t(z,w)|$ denote the length of
$\lambda_t(z,w)$, i.e. the distance between $z$ and $w$ in the flat
metric on $\tilde{\gamma}(t)$. Suppose $\lambda_0(z,w)$ is not an edge
in $W$. Then, either $\lambda_0(z,w)$ is a saddle connection not in
$W$, in which case $|\lambda_0(z,w)| \ge \sqrt{2} \ell_{min}(q)$, or
$\lambda_0(z,w)$ is a union of at least two saddle connections, so
that $|\lambda_0(z,w)| \ge 2 \ell_{min}(q) \ge \sqrt{2}\ell_{min}(q)$. 
It now follows
from (\ref{eq:change:in:length}) that for all $t \in [0,\rho_1]$, 
\begin{displaymath}
|\lambda_t(z,w)| \ge (\sqrt{2} - n m t) \ell_{min}(q),
\end{displaymath}
where $n$ is the maximum intersection number of $\Delta'$ with a
saddle connection in the Delaunay triangulation of $\tilde{S}$, 
and $m$ is the maximal number of saddle
connections in $\lambda_t(z,w)$. Note that both $n$ and $m$ are
bounded by the genus. Now we choose $\rho_1$ so that $\sqrt{2} - n m \rho_1 \ge
(1+\rho_1)^2$. Then, if $\lambda_0(z,w)$ is not a saddle connection 
in $W$, then for all $t \in [0,\rho_1]$,
\begin{equation}
\label{eq:length:other:curves}
|\lambda_t(z,w)| \ge \ell_{min}(\gamma(0)) \sqrt{1+t^2}. 
\end{equation}
Now (\ref{eq:shortest:increases}) follows from (\ref{eq:length:UcapW})
and (\ref{eq:length:other:curves}). 

We now estimate the Hodge length of the path $\gamma: [0,\rho_1] \to
\qonetg$. By Lemma~\ref{lemma:flat:to:hyperbolic} and
Lemma~\ref{lemma:no:flat:cylinders}, 
\begin{displaymath}
\ell_{\gamma(t)}(\sigma) \ge \frac{C_g} {|\log( \ell_{min}(q) \sqrt{1+t^2})|},
\end{displaymath}
where $\ell_\sigma(S)$ is defined to be the infimum over all simple
closed curves $\alpha$ of $\ell_{\alpha}(\sigma)$ (and $\sigma$ is the
hyperbolic metric in the conformal class of $S$). 
By construction, the intersection number  of $\Delta'$ with a short
basis (see \S\ref{sec:modified:norm}) 
of $\gamma(t)$ is bounded depending only on the genus.
Therefore, by (\ref{eq:expand:in:basis}) and
(\ref{eq:def:modified:hodge:norm}), 
\begin{displaymath}
\|\gamma'(t)\|'  \le {C}{|\log( \ell_{min}(q) \sqrt{1+t^2})|^{1/2}},
\end{displaymath}
where $C$ depends only on the genus. Thus, if $q_1 =
\gamma(\rho_1)$, then
\begin{displaymath}
d_H(q,q_1) \le \int_0^{\rho_1}   \|\gamma'(t)\|' \,dt \le
\int_{0}^{\rho_1} {C}{|\log( \ell_{min}(q) \sqrt{1+t^2})|^{1/2}} \, dt.
\end{displaymath}
Thus (\ref{eq:hodge:distance:bound}) holds. 
\qed\medskip

\begin{lemma}
\label{lemma:opening:up}
There exists $\epsilon_0 > 0$ (depending on $\cK$) 
such that for all $\epsilon < \epsilon_0$ and for all $q \in \pi^{-1}(\cK)$
with $\ell_{min}(q) < \epsilon$, there exists $q'$ on the same leaf of
$\cF^{uu}$ as $q$ with $\ell_{min}(q') \ge \epsilon$, and $d_H(q,q') \le k
\epsilon |\log \epsilon|^{1/2}$, where $k$ depends only on $\cK$. 
\end{lemma}

\bold{Proof.} We define a sequence $q_n$ as follows: let $q_0 = q$.
If $q_n$ has been defined already, we apply
Lemma~\ref{lemma:expand:one:step} with $q = q_n$, and define $q_{n+1}$
to be the point $q'$ guaranteed by
 Lemma~\ref{lemma:expand:one:step}. We obtain a
sequence $q_n$ with 
\begin{displaymath}
\ell_{min}(q_{n+1}) \ge \ell_{min}(q_n)(1+\rho_1^2)^{1/2} = \rho_2 \ell_{min}(q_n),
\end{displaymath}
where we let $\rho_2 = (1+\rho_1^2)^{1/2}$. 
Define $t_n$ inductively by $t_0 = 0$, $t_{n+1} = t_n + \rho_1
\ell_{min}(q_n)$, $n \ge 0$. Then,
\begin{displaymath}
t_n = \rho_1 \sum_{k=0}^{n-1} \ell_{min}(q_k) \le \rho_1 \sum_{k=0}^{n-1}
\frac{\ell_{min}(q_n)}{ \rho_2^{n-k}} = \rho_1 \ell_{min}(q_n) \frac{\rho_2^{-1}
- \rho_2^{-n-1}}{1 -\rho_2^{-1}} \le \frac{ \rho_1 \ell_{min}(q_n)}{\rho_2 -1}.
\end{displaymath}
Then, 
\begin{displaymath}
t_{n+1} =t_n + \rho_1 \ell_{min}(q_n) \le  \frac{\rho_1 \rho_2}{\rho_2 - 1} \ell_{min}(q_n).
\end{displaymath}
Thus, for $t \in [t_n,t_{n+1}]$,
\begin{displaymath}
\ell_{min}(q_n) \ge \frac{\rho_2 -1}{\rho_1 \rho_2} t \equiv \rho_3 t.
\end{displaymath}
We have 
\begin{displaymath}
d_H(q_n,q_{n+1}) \le \int_{t_n}^{t_{n+1}}  {C}{|\log( \ell_{min}(q_n)
  \sqrt{1+t^2})|^{1/2}} \, dt \le  \int_{t_n}^{t_{n+1}} {C}{|\log(
  \rho_3 t \sqrt{1+t^2})|^{1/2}} \, dt.
\end{displaymath}
Thus, 
\begin{displaymath}
d_H(q_0,q_n) \le \int_{0}^{t_n} {C}{|\log( \rho_3 t
  \sqrt{1+t^2})|^{1/2}} \, dt.
\end{displaymath}
We choose $n$ so that $\ell_{min}(q_n)$ is comparable to $\epsilon$. 
Let $q' =\gamma(t_n)$. Then, $\ell_{min}(q') \ge t_n \approx
\epsilon$. Now the modified Hodge length of the path $q_0,q_1, \dots
,q_n$ is $O(\epsilon |\log \epsilon|^{1/2})$. 
\qed\medskip

\bold{Proof of Theorem~\ref{theorem:hodge:equivalent:to:euclidean}.}
Since $\cK$ is compact, the intersection of the multiple zero locus
with $\pi^{-1}(\cK)$ is (contained in) finite union of hyperplanes $H_1, \dots,
H_n$. Each hyperplane $H_j$ has complex codimension $1$. We can choose
$\epsilon_0 > 0$ such that any two $H_j$ which do not intersect in
$\pi^{-1}(\cK)$ are at least $\epsilon_0$ apart. Let $\delta_0 > 0$ be
a lower bound on the angle between any two  $H_j$ which do intersect
in $\pi^{-1}(\cK)$. Clearly $\epsilon_0$ and $\delta_0$ depend only on $\cK$. 

Let $Z$ be the locus where our quadratic differential has either a
zero of order at least $3$ 
or two zeroes each of order at least $2$. Then there exists a constant
$k_0$ (depending only on $\delta_0$ and thus only on $\cK$) such that
for any $q \in \pi^{-1}(\cK)$ and any $\epsilon > 0$, if the ball
$B_E(q,\epsilon)$ intersects at least two hyperplanes $H_j$ then
$d_E(q,Z) < k_0 \epsilon$. 


Take two points $q_1, q_2 \in \pi^{-1}(\cK)$ on the same leaf of
$\cF^{uu}$, with 
$d(q_1,q_2) = \epsilon$. Choose $k_3 > k_0$ (thus $k_3$ depends only on
$\cK$). We also assume that $k_3 \epsilon < \epsilon_0$. 
By Lemma~\ref{lemma:opening:up} (with $k_3 \epsilon$ in
place of $\epsilon$), there exist $q_1', q_2'$ so that for $i=1,2$, 
we have $d_H(q_i,q_i') \le k_1 \epsilon |\log
\epsilon|^{1/2}$, $d_E(q_i,q_i') \le k_2 \epsilon$, and $\ell_{min}(q_i) = k_3
\epsilon$, where $k_1, k_2$ depend only on $\cK$. Since $k_3 > k_0$
and $k_3 \epsilon < \epsilon_0$, 
there exist $\epsilon_1$, $\epsilon_1'$ such that
$k_3 \epsilon/2 > \epsilon_1 > \epsilon_1' > \epsilon/(2k_3)$, and
for all $q \in \pi^{-1}(\cK)$, either $d_E(q,Z) < \epsilon_1$ or
$B_E(q,\epsilon_1')$ contains at most one hyperplane from the multiple
zero locus. 

Note that the intersection of
$Z$ with $\pi^{-1}(\cK)$ has complex codimension at least $2$. Hence, the
intersection of the $\epsilon_1$-neighborhood of $Z$ with 
$\cF^{uu}(q_1) \cap \pi^{-1}(\cK)$ 
is contained in the $\epsilon_1$-neighborhood of 
a finite union of hyperplanes, each of real codimension at least $2$. 
Then there exists a constant $k_4$ depending only on $\cK$, and 
a path $\gamma$ connecting $q_1'$ to $q_2'$ 
of length at most $k_4 \epsilon$, which avoids the
$\epsilon_1$-neighborhood of $Z$.

Now let $p_0 = q_1'$, and mark points $p_i$ along $\gamma$ which 
are $\epsilon_1'/2$ apart in the Euclidean metric. 
We have $p_n = q_2'$. Let $B_i$ be the 
ball of Euclidean diameter $\epsilon_1'/2$ which contains $p_i$ and $p_{i+1}$
on its boundary. By construction $B_i$ contains at most one hyperplane
(which we will denote $L$) from the multiple zero locus. 
Note that by Lemma~\ref{lemma:flat:to:hyperbolic}, 
Lemma~\ref{lemma:no:flat:cylinders} and (\ref{eq:star:c:beta:prime:norm}), 
for any $p \in B_i$ and any tangent vector
$v$ at $p$, the modified Hodge norm of $v$ can be estimated as
\begin{equation}
\label{eq:crude:estimate:of:hodge:norm}
\|v\|_H' \le C |\log d_E(p,L)|^{1/2} \|v\|_E,
\end{equation}
where $\|v\|_E$ is the Euclidean norm of $v$, and $d(p,L)$ denotes the
Euclidean distance between the point $p$ and the hyperplane $L$. 

Let $p_i'$ be the farthest point in
$B_i$ from the hyperplane. Then, after connecting $p_i$ and $p_i'$ by
a straight line path and using
(\ref{eq:crude:estimate:of:hodge:norm}), we see that $d_H(p_i,p_i') = O(\epsilon |\log
\epsilon|^{1/2})$, and also $d_H(p_i',p_{i+1}) = O(\epsilon |\log
\epsilon|^{1/2})$. Thus, since the number of $B_i$ along the path is
bounded by a constant depending only on $\cK$, we finally obtain
\begin{displaymath}
d_H(q_1',q_2') = O(d_E(q_1',q_2') |\log \epsilon|^{1/2}).
\end{displaymath}
\qed\medskip

\subsection{The Non-expansion and Decay of the Hodge Distance.}
\label{sec:hodge:distance}

\begin{theorem}
\label{theorem:hodge:distance:nonincrease}
Suppose $q \in \qonetg$ and $q' \in \qonetg$ 
are in the same leaf of $\cF^{ss}$. Then
\begin{itemize}
\item[{\rm (a)}] There exists a constant $c_H > 0$ such that for all
  $t \ge 0$, 
\begin{displaymath} d_H(g_t q, g_t q') \le c_H d_H(q,q').
\end{displaymath}
\item[{\rm (b)}] Suppose $\epsilon > 2 \epsilon_1$, where $\epsilon_1$
  is as in the definition of short basis (see \S\ref{sec:modified:norm}). Let $K = \{ q \in \qonetg
  \st \ell_{min}(q) > \epsilon \}$. Suppose $d_H(q,q') < 1$, and $t > 0$ 
is such that 
\begin{equation}
\label{eq:at:least:half:in:K}
|\{ s \in [0,t] \st g_s q \in \Gamma K \}| \ge (1/2) t.
\end{equation}
Then for all $0 < s < t$, 
\begin{displaymath}
d_H(g_s q, g_s q') \le C e^{-c s} d_H(q,q'),
\end{displaymath}
where $c$ and $C$ depend only on $g$ and $\epsilon$. 
\end{itemize}
\end{theorem}

\bold{Proof.} The statement  (a) follows immediately from
Theorem~\ref{theorem:nonincrease:modified:norm}. For the second
statement, let $\gamma: [0,\rho] \to \qonetg$ be a modified
Hodge length minimizing path connecting $q$ to $q'$ (and staying in
the same leaf of $\cF^{ss}$). We assume that $\gamma$ is parametrized so that
$d_H(\gamma(u),q) = u$. By assumption, $\rho < 1$. 

Let $N= 4 c_H/(c_1 \epsilon\rho)$, where $c_1$
is as in Theorem~\ref{theorem:hodge:equivalent:to:euclidean}. 
For $0 \le j \le N$, let $u_j = j c_1 \epsilon/(4 c_H)$. We now claim
that for all $j$ there exists $t_j$ such that for $s > t_j$ such that
$g_s q \in \Gamma K$, we have for all $0 \le u < u_j$, 
\begin{equation}
\label{eq:ind:euc:bound}
d_E(g_s \gamma(u), g_s q ) < \epsilon \sum_{k=0}^j 2^{-k}, 
\end{equation}
and
\begin{equation}
\label{eq:ind:hodge:bound}
d_H(g_s \gamma(u), g_s q) < C e^{-c s} u, 
\end{equation}
where $C$ and $c$ depend only on $g$ and $\epsilon$. 
The equations (\ref{eq:ind:euc:bound}) and (\ref{eq:ind:hodge:bound})
will be proved by induction on $j$. If $j=0$ then there is nothing to 
prove. Now assume $j \ge 1$, and (\ref{eq:ind:euc:bound}) and
(\ref{eq:ind:hodge:bound}) are true for all $u \le u_{j-1}$. 
Suppose $u < u_j$. 
By (a), we have for all $t > 0$, 
\begin{displaymath}
d_H(g_t \gamma(u), g_t
\gamma(u_{j-1})) < c_H (u-u_{j-1}) < c_1 \epsilon/4.
\end{displaymath}
Therefore by Theorem~\ref{theorem:hodge:equivalent:to:euclidean}, $d_E(g_t
\gamma(u), g_t \gamma(u_{j-1})) < \epsilon/2$. 
Therefore, by the inductive assumption (\ref{eq:ind:euc:bound}), 
for the $s > t_{j-1}$ such that $g_s q \in \Gamma K$, 
for all $u_{j-1} < u < u_j$, we have $\ell_{min}(g_s \gamma(u))
> \epsilon/2$. Thus, by (\ref{eq:modified:hodge:compare:to:hodge}), 
for all $u_{j-1} < u < u_j$, the modified Hodge
norm of $\gamma'(u)$ is within a constant of the Hodge norm of
$\gamma'(u)$. Then by Theorem~\ref{theorem:forni}, 
for all $u_{j-1} < u < u_j$, 
\begin{equation}
\label{eq:tmp:exp:decay}
d_H(g_t \gamma(u), g_t \gamma(u_{j-1})) < C e^{-c t} d_H(\gamma(u),
\gamma(u_{j-1})). 
\end{equation}
Therefore, by
Theorem~\ref{theorem:hodge:equivalent:to:euclidean}, there exists $t_j
> 0$ (depending only on $g$ and $\epsilon$)  
such that for $s > t_j$ with $g_s q \in
\Gamma K$, and $u_{j-1} \le u < u_j$, (\ref{eq:ind:euc:bound}) holds. Also
(\ref{eq:ind:hodge:bound}) follows from (\ref{eq:tmp:exp:decay}). 
The induction terminates after finitely many steps (depending only on
$g$ and $\epsilon$). Thus (b) holds. 
\qed\medskip

\section{The multiple zero locus.}
\label{sec:multiple:zeroes}
In this section we prove Theorem~\ref{theorem:mz:volume} and
Theorem~\ref{theorem:badgeodesic:volume}. 

Let $\cK$ be a compact subset of $\cT_g$. We assume that $\cK$ is
contained in one fundamental domain for the action of $\Gamma$
on $\cT_g$ (and thus we can identify $\cK$ with a subset of 
$\cM_g$). 
All of our implied constants will depend on $\cK$. 

\bold{Notation.} 
Suppose $W \subset \qonetg$ and $s > 0$. 
Let $W(s)$ denote the set of $q \in \qonetg$ such that
there exists $q' \in W$ on the same leaf of $\cF^{uu}$ as $q$ such 
that $d_H(q,q') < s$. Recall that for a subset $A \subset \cT_g$, 
$Nbhd_{r}(A)$ denotes the set of points within Teichm\"uller
distance $r$ of $A$. We will also use $A(r)$ to denote $Nbhd_{r}(A)$.
Let $\bar{\nu}$ and $\eta^+$ be as in \S\ref{sec:conditional}, so that
so that for a set $F$ contained in a leaf of $\cF^{uu}$, 
$\bar{\nu}(\eta^+(F)) = \mu_{\alpha^{uu}}(F)$. 

\begin{lemma}
\label{lemma:technical:projection}
Suppose $U \subset \pi^{-1}(\cK)$, $\delta > 0$ and $t > 0$. Let $W = g_t U \cap
\Gamma \pi^{-1}(\cK)$. There is a $C(\delta)>0$ so that
\begin{equation}
\label{eq:tech:proj:one}
\mm(Nbhd_2(\pi(W))) \le C(\delta)
\overline{\nu}(\eta^+(W(\delta))), 
\end{equation}
(By our convention, the constant $C(\delta)$ depends on $\cK$ as well). Also
\begin{equation}
\label{eq:tech:proj:two}
\overline{\nu}(\eta^+(W(1))) \le C(\delta) \overline{\nu}(\eta^+(W(\delta))).
\end{equation}
\end{lemma}

\bold{Proof.} Let $\delta_0 = \delta_0(\cK,\delta)$ be a constant to be
chosen later. 
We decompose $U$ into pieces $U_\alpha$ such that each piece is within
(modified) Hodge distance $\delta_0/2$ of a single leaf of
$\cF^{uu}$. Let
\begin{displaymath}
W_\alpha = g_t U_\alpha \cap \Gamma \pi^{-1}(\cK). 
\end{displaymath}
Then
\begin{displaymath}
\mm(Nbhd_2(\pi(W))) \le \sum_\alpha \mm(Nbhd_2(\pi(W_\alpha)))
\end{displaymath}
In view of Theorem~\ref{theorem:hodge:equivalent:to:euclidean} the
number of pieces is bounded depending only on $\cK$ and $\delta$
(since $\delta_0 = \delta_0(\cK,\delta)$). Thus, it is enough to show that
(\ref{eq:tech:proj:one}) holds with $W$ replaced by $W_\alpha$. 

We now claim that without loss of 
generality, we may assume that $U_\alpha$ has the
following ``product property'': 
given $q_1$, $q_2 \in U_\alpha$, there exists $q_2'
\in U_\alpha$ on the same leaf of $\cF^{uu}$ as $q_1$ and on the same
leaf of $\cF^{s}$ as $q_2$. If not, 
let 
\begin{displaymath}
U_\alpha' = U_\alpha \cup \{ \cF^{uu}(q_1) \cap \cF^{s}(q_2) \st
\text{ $q_1, q_2 \in 
U_\alpha$} \}. 
\end{displaymath}
Then $\eta^+(U_\alpha) = \eta^+(U_\alpha')$,
and therefore $\eta^+(g_t U_\alpha) = \eta^+(g_t U_\alpha')$.
Also if $\delta'=\delta'(\delta, K)$ is sufficiently small then
$W_\alpha' = g_t U_\alpha' \cap \Gamma \pi^{-1}(\cK)$ satisfies
$W_\alpha'(\delta') \subset W_\alpha(\delta)$. Therefore we can proceed
with the rest of the proof with $\delta'$ instead of $\delta$ and
$U_\alpha'$ instead of $U_\alpha$. Therefore, without loss of
generality, we may assume that $U_\alpha$ has the product property. 
Therefore, $g_t U_\alpha$ also has the
product property. 

Pick a maximal $\Delta \subset \pi(W_\alpha)$ such that for any
two distinct $X,Y \in \Delta$, $d_{\cT}(X,Y) = 1$. Then,
\begin{displaymath}
Nbhd_2(\pi(W_\alpha)) \subset \bigcup_{X \in \Delta} B_{\tau}(X,3), 
\end{displaymath}
and hence,
\begin{equation}
\label{eq:m:Nbhd2:Delta}
\mm(Nbhd_2(\pi(W_\alpha))) \le \mm \left( \bigcup_{X \in \Delta}
  B_{\tau}(X,3) \right) \le \sum_{X \in \Delta} \mm(B(X,3)) \le C(\cK)
|\Delta|,
\end{equation}
where $|\Delta|$ denotes the cardinality of $\Delta$, and we have used
the fact that $\Delta \subset \Gamma\cK$.

For each $X \in \Delta$, pick one $q \in W_\alpha \cap \pi^{-1}(X)$. Let
$\Delta' \subset W_\alpha$ denote the resulting set of $q$'s. 
Let
$B_E^{uu}(q,r)$ denote the set of $q' \in \qonetg$ on the same leaf of
$\cF^{uu}$ as $q$ with $d_E(q,q') < r$. 
We claim
that for $\delta_0$ sufficiently small, we can pick $\delta_2$
depending only on $\cK$ such that for all distinct
pairs $q_1, q_2 \in \Delta'$, 
\begin{equation}
\label{eq:eta:plus:Buu:disjoint}
\eta^+(B_E^{uu}(q_1,\delta_2)) \cap \eta^+(B_E^{uu}(q_2,\delta_2)) = \emptyset.
\end{equation}
To prove (\ref{eq:eta:plus:Buu:disjoint}), 
suppose $q_1,
q_2 \in \Delta'$, and $q_1 \ne q_2$. Let $q_2'$ be such that $q_1$ and
$q_2'$ are on the same leaf of $\cF^{uu}$ and $q_2$ and $q_2'$ are on
the same leaf of $\cF^{s}$. By
Theorem~\ref{theorem:hodge:distance:nonincrease}, we have
$d_H(q_2,q_2') < c_H \delta_0$. Therefore, by
Theorem~\ref{theorem:hodge:equivalent:to:euclidean} and
Lemma~\ref{lemma:Teich:equivalent:to:euc}, we can choose $\delta_0$
small enough so that $d_\cT(\pi(q_2), \pi(q_2')) \le 1/5$. Hence, 
$\pi(q_2') \subset
\Gamma \cK'$ where $\cK' \subset \cT_g$ is compact, 
and by the triangle inequality,
\begin{equation}
\label{eq:tmp:triangle}
d_{\cT}(\pi(q_1), \pi(q_2')) \ge d_{\cT}(\pi(q_1), \pi(q_2))  -
d_{\cT}(\pi(q_2), \pi(q_2'))  \ge 4/5. 
\end{equation}
By Lemma~\ref{lemma:Teich:equivalent:to:euc}, there exists $\delta_1 >
0$ (depending only on $\cK'$) 
such that for all $q$, $q' \in \cK'$ on the same leaf of $\cF^{uu}$
with $d_E(q,q') < \delta_1$, we have $d_{\cT}(\pi(q), \pi(q')) < 1/5$.
Then, by (\ref{eq:tmp:triangle}) and the triangle inequality,
$B_E^{uu}(q_1,\delta_1)$ and 
$B_E^{uu}(q_2',\delta_1)$ are disjoint. It follows that as subsets of
$\pmf$, 
\begin{equation}
\label{eq:pmf:disjoint}
\eta^+(B_E^{uu}(q_1,\delta_1)) \cap \eta^+(B_E^{uu}(q_2',\delta_1)) = \emptyset
\end{equation}
Let $\cF^{s}(q)$ denote the leaf of $\cF^{s}$ through $q$. Since
$\eta^+$ is continuous and $\cK'$ is compact, 
there exist constants $0 < \delta_2 < \delta_1$ and
$\delta_3 > 0$ (depending only on $\cK'$) such that for all $q,q' \in
\cK'$ with $q \in \cF^s(q')$ and $d_H(q,q') < \delta_3$, we have, as
subsets of $\pmf$, 
\begin{displaymath}
\eta^+(B_E^{uu}(q,\delta_2)) \subset \eta^+(B_E^{uu}(q', \delta_1)). 
\end{displaymath}
We now choose $\delta_0 < \delta_3/c_H$. Then, since $q_2$ and $q_2'$ are on
the same leaf of $\cF^s$ and $d_H(q_2,q_2') < c_H \delta_0 < \delta_3$, we have
\begin{displaymath}
\eta^+(B_E^{uu}(q_2,\delta_2)) \subset \eta^+(B_E^{uu}(q_2', \delta_1)). 
\end{displaymath}
Now using
(\ref{eq:pmf:disjoint}), we get that as subsets of
$\pmf$,
\begin{displaymath}
\eta^+(B_E^{uu}(q_1,\delta_2)) \cap \eta^+(B_E^{uu}(q_2,\delta_2))
\subset \eta^+(B_E^{uu}(q_1,\delta_1)) \cap
\eta^+(B_E^{uu}(q_2',\delta_1)) = \emptyset.
\end{displaymath}
This completes the proof of (\ref{eq:eta:plus:Buu:disjoint}). 

By Theorem~\ref{theorem:hodge:equivalent:to:euclidean}, there exists
$0 < \delta_3 < \delta_2$ such that for all $q \in \Delta'$ and all $q' \in
\cF^{uu}(q)$ with $d_E(q,q') < \delta_3$ we have $d_H(q,q') <
\delta$. For $q \in \Delta'$, let 
\begin{displaymath}
H(q) = \eta^+(B_E(q,\delta_3)) \subset \eta^+(W(\delta)),
\end{displaymath}
Consider the collection of
``balls'' $\{ H(q) \st q \in \Delta' \}$. By
(\ref{eq:eta:plus:Buu:disjoint}) the sets $H(q)$ are pairwise
disjoint viewed as subsets of $\pmf$ (or alternatively the subsets
$Cone(H(q)) \subset \mf$ intersect only at the origin). Also, by the
definition of the Thurston measure $\nu$ and the compactness of $\cK$, 
there exists a constant $c = c(\cK,\delta)$ such that
\begin{displaymath}
\bar{\nu}(H(q)) \ge c, \quad\text{ for all $q \in \Delta'$.}
\end{displaymath}
Hence, 
\begin{equation}
\label{eq:bar:nu:W:delta:Delta}
\overline{\nu}(\eta^+(W(\delta))) \ge \overline{\nu}\left(\bigcup_{q \in
    \Delta'} H(q) \right)  = 
\sum_{q \in \Delta'}
\overline{\nu}(H(q)) \ge c |\Delta'|  = c |\Delta|.
\end{equation}
Now (\ref{eq:tech:proj:one}) (with $W$ replaced by $W_\alpha$) 
follows from (\ref{eq:m:Nbhd2:Delta}) and
(\ref{eq:bar:nu:W:delta:Delta}). 
Finally, 
$$\overline{\nu}(\eta^+(W(1))) \le C(\delta,\cK)
\overline{\nu}(\eta^+(W(\delta))),$$ since $d_H$ is equivalent to
$d_E$. 
\qed\medskip

\bold{The sets $K_i$ and $U_i$.}
Let $K_1 \subset \qonemg$ be a compact set. (In our application, 
$K_1$ will be chosen disjoint from the multiple zero locus). 
Let $K_3 \subset K_2 \subset K_1$ and $1 \ge \delta > 0$ be
such that if $q \in K_i$ and $d_H(q,q') < c_H \delta$ then $q' \in
K_{i-1}$, where $c_H$ is as in
Theorem~\ref{theorem:hodge:distance:nonincrease} (a). 
We assume $\overline{\mu}(K_3) > (1/2)$, 
where $\overline{\mu}$ is the normalized Lebesgue measure on
$\qonemg$.
For $T_0 > 0$ let $U_i = U_i(T_0)$ be the set of $q \in \qonemg$ 
such that there exists  $T > T_0$ so that
\begin{displaymath}
|\{ t \in [0, T] \st g_t q \in K_i^c \}| \geq (1/2) T.
\end{displaymath}
Then, for all $T > T_0$ and all $q \not\in U_i$, 
\begin{equation}
\label{eq:outside:Ki}
|\{ t \in [0,T] \st g_t q \in K_i^c \}| < (1/2)T.
\end{equation}
From the definition, we have $U_1 \subset U_2 \subset U_3$. 
By the ergodicity of the geodesic flow, for every $\theta > 0$ there exists
$T_0 > 0$ such that $\overline{\mu}(U_3) < \theta$.

Let $\cK_1 = \cK$. We can choose compact subsets $\cK_0$, $\cK_2$,
$\cK_3$ of $\cT_g$ such that for $0 \le i \le 2$, any $q \in
\pi^{-1}(\cK_i)$ and any $q'$ on the same leaf of $\cF^{uu}$ as $q$ 
with $d_H(q,q') < c_H$, we have $q' \in \pi^{-1}(\cK_{i+1})$. 
By (the lower bound in) Theorem~\ref{theorem:hodge:equivalent:to:euclidean}
and Lemma~\ref{lemma:Teich:equivalent:to:euc}, each $\cK_i$ is compact. 
Let $U_i' = p^{-1}(U_i) \cap \pi^{-1}(\cK_i)$,
where $p$ is the natural map from $\qonetg$ to $\qonemg$.
Note that $U_1' \subset U_2' \subset U_3' \subset \qonetg$.

\begin{lemma}
\label{lemma:small:measure:for:bad:set}
Let $U_i'$, $1 \le i \le 3$ 
be as in the above paragraph. Then, for all $t > 0$, 
\begin{displaymath}
\mm(Nbhd_2(\pi(g_t U_1') \cap \Gamma \cK)) \le C(\delta) e^{ht}
\overline{\nu}(\eta^+(U_2')), 
\end{displaymath}
and if $W(1)$ is defined as in
Lemma~\ref{lemma:technical:projection} with $W = g_t U_1' \cap \pi^{-1}(\Gamma
\cK)$, then 
\begin{displaymath}
\overline{\nu}(\eta^+(W(1))) \le C'(\delta) e^{ht}
\overline{\nu}(\eta^+(U_2')). 
\end{displaymath}
In particular (c.f. Theorem~\ref{theorem:sx:multiple:zero:locus}), 
for any $\epsilon > 0$ is is possible to choose
$T_0$ such that if $U_1$ is defined by (\ref{eq:outside:Ki}) and
$U_1' = U_1 \cap \pi^{-1}(\cK_1)$ then for all $t > T_0$,
\begin{equation}
\label{eq:m:Nbhd:delta:U1:small:measure}
\mm(Nbhd_2(\pi(g_t U_1')) \cap \Gamma \cK) \le \epsilon e^{ht},
\end{equation}
and for $W = g_t U_1' \cap \pi^{-1}(\Gamma \cK)$ we have
\begin{equation}
\label{eq:Wone:small:measure}
\overline{\nu}(\eta^+(W(1))) \le \epsilon e^{ht}.
\end{equation}
\end{lemma}

\bold{Proof.} We will apply Lemma~\ref{lemma:technical:projection}
to the set $W = g_t U_1' \cap \pi^{-1}(\Gamma \cK)$. 
We claim that $W(\delta) \subset g_t
U_2'$. Indeed, suppose  $g_t q' \in W(\delta)$. Then 
there exists $g_t q  \in W$ with
\begin{displaymath}
d_H(g_t q, g_t q') < \delta. 
\end{displaymath}
Then, by Theorem~\ref{theorem:hodge:distance:nonincrease} (a), 
for all $0 \le s \le t$,
\begin{equation}
\label{eq:tmp:nonexp}
d_H(g_s q, g_s q') < c_H \delta
\end{equation}
Since $W = g_t U_1'$, $g_t q \in W$ implies $q \in U_1'$. 
Then, assuming $t > T_0$, 
for at least half the values of $s \in [0,t]$, 
\begin{equation}
\label{eq:tmp2:gsq:not:in:K}
g_s q \not\in K_1.
\end{equation}
Then, (\ref{eq:tmp:nonexp}) and the definition of $K_2$ 
imply that for the $s \in [0,t]$ for which
(\ref{eq:tmp2:gsq:not:in:K}) holds, 
$g_s q' \not\in K_2$. 
This implies $q' \in U_2$. 
Since $d_H(q,q') < c_H \delta < c_H$ and $q \in \pi^{-1}(\cK) =
\pi^{-1}(\cK_1)$, 
we have $q' \in \pi^{-1}(\cK_2)$. Thus, $q' \in U_2 \cap
\pi^{-1}(\cK_2) = U_2'$, and so $g_t q' \in g_t U_2$.
This implies the claim, and thus the first two statements
of the Lemma. 

The same argument as the proof of the claim shows that
if $q \in U_2'$ and $q' \in \cF^{ss}(q) \cap \pi^{-1}(\cK)$ 
with $d_H(q,q') < \delta$ then $q' \in U_3'$. This (together with
Theorem~\ref{theorem:hodge:equivalent:to:euclidean}) implies
that there exists $C_1(\delta)$ such that \mc{What does $C_1(\delta)$ depend on?}
\begin{displaymath}
\overline{\nu}(\eta^+(U_2')) \le C_1(\delta) \overline{\mu}(U_3').
\end{displaymath}
Hence if we choose $T_0 > 0$ so that $\overline{\mu}(U_3) < C(\delta)
C_1(\delta) \epsilon$, then (\ref{eq:m:Nbhd:delta:U1:small:measure}) 
follows. Similarly, if we choose $T_0 > 0$ so that in addition 
$\overline{\mu}(U_3) < C'(\delta) C_1(\delta) \epsilon$, 
then (\ref{eq:Wone:small:measure}) follows. 
\qed\medskip

\bold{Proof of Theorem~\ref{theorem:badgeodesic:volume}.}
Let $K' = K_1$, and let $T_0$, $U_1$ and $U_1'$ be as in
Lemma~\ref{lemma:small:measure:for:bad:set}. Then, for $R > T_0$,
and $X\in \cK$,
\begin{displaymath}
B_R(X,\cK,K') \subset \bigcup_{0 \le t \le R} \pi(g_t U_1')\cap
\Gamma \cK \quad \subset \quad \bigcup_{n=0}^{\lfloor
  R \rfloor} \bigcup_{t \in [n,n+1]} \pi(g_t U_1') 
\cap \Gamma \cK,
\end{displaymath}
where $\lfloor x \rfloor$ denotes the integer part of $x$.
Then,
\begin{displaymath}
\mm(Nbhd_1(B_R(X,\cK,K'))) \le \sum_{n=0}^{\lfloor R \rfloor}
\mm(Nbhd_2(\pi(g_n U_1')\cap \Gamma \cK)) 
\le C \epsilon \sum_{n=0}^{\lfloor R \rfloor}  e^{hn},
\end{displaymath}
where we have used (\ref{eq:m:Nbhd:delta:U1:small:measure}).
Since $\epsilon$ is arbitrary,
Theorem~\ref{theorem:badgeodesic:volume} follows.
\qed\medskip

\begin{lemma}
\label{lemma:small:measure:for:nbhd:of:multiple:zeroes}
Suppose $V \subset \pi^{-1}(\cK)$ and $\delta' > 0$. 
Then, for $t$ sufficiently large (depending on $\cK$, $V$ 
and $\delta'$),
\begin{displaymath}
\mm(Nbhd_{1}(\pi(g_t V) \cap \Gamma \cK_0)) 
\le C \overline{\nu}(\eta^+(V(\delta')))e^{ht},
\end{displaymath}
where $C$ depends only on $\cK$.
\end{lemma}

\bold{Proof.} 
Let $Y = g_t V \cap \Gamma \cK_0$. Choose $T_0$ so that
(\ref{eq:Wone:small:measure}) holds for $t > T_0$ and
$\overline{\nu}(\eta^+(V(\delta')))$ instead of $\epsilon$. 
Let $U_1$, $U_1'$ and $W = g_t U_1' \cap \pi^{-1}(\cK)$ be as in
Lemma~\ref{lemma:small:measure:for:bad:set}, so in particular, for $t
> T_0$, 
\begin{equation}
\label{eq:W1new:small:measure}
\overline{\nu}(\eta^+(W(1))) \le \overline{\nu}(\eta^+(V(\delta'))) 
e^{ht}.
\end{equation}
We claim that there exist $T_1 > T_0$, 
depending only on $\cK$ such
that for $t > T_1$, 
\begin{equation}
\label{eq:tmp:claim}
Y(1) = W(1) \cup g_t (V(\delta')). 
\end{equation}

Indeed, if $g_t q \in Y(1)$ then by definition there exists $q' \in V
\subset \cK$ with $d_H(g_t q, g_t q') < 1$, and $q'$ is on the same
leaf of $\cF^{uu}$ as $q$. 
We consider two cases: either $q' \in U_1$
or $q' \not\in U_1$. If $q' \in U_1$ then $q \in U_1' = U_1 \cap
\cK$, hence $g_t q' \in g_t U_1' \cap \Gamma \cK = W$. Hence
in this case, $g_t q \in W(1)$. 
If $q' \not\in U_1$, then by (\ref{eq:outside:Ki}) 
and Theorem~\ref{theorem:hodge:distance:nonincrease} (b), 
\begin{displaymath}
d_H(q,q') = d_H(g_{-t} g_t q, g_{-t} g_t q') \le C e^{-ct} d_H(g_t q,
g_t q') \le C e^{-ct}.
\end{displaymath}
We choose $T_1 > T_0$ so that $C e^{-c T_1} < \delta'$. 
Thus, in this case, for $t > T_1$, $q \in V(\delta')$ and hence
(\ref{eq:tmp:claim}) follows. 

Now, for $t > T_1$, 
\begin{displaymath}
\overline{\nu}(\eta^+(Y(1))) \le \overline{\nu}(\eta^+(W(1))) +
\overline{\nu}(\eta^+(g_t(V(\delta')))) 
 \le 2 \overline{\nu}(\eta^+(V(\delta'))) e^{ht}
\end{displaymath}
where we have used (\ref{eq:W1new:small:measure}). We now apply
Lemma~\ref{lemma:technical:projection} with $\delta = 1$. 
The lemma follows from (\ref{eq:tech:proj:one}). 
\qed\medskip

\bold{Proof of Theorem~\ref{theorem:mz:volume}.}
This follows immediately
from Lemma~\ref{lemma:small:measure:for:nbhd:of:multiple:zeroes},
and the fact that we can choose an relatively open 
$V \subset \cK$ and $\delta' > 0$ such that
$V$ contains the intersection of $\cK$ with
the multiple zero locus, and $\overline{\nu}(\eta^+(V(\delta')))$ is
arbitrarily small, (see Theorem~\ref{theorem:sx:multiple:zero:locus}).

\section{Volume Asymptotics.}
\label{sec:volume}
In this section we prove Theorem~\ref{theorem:volume:growth}.
Let $F_R(X,Y) = |\Gamma \cdot Y \cap B_R(X)|$ be as in  \S\ref{sec:counting}. 
We need the following:
\begin{theorem}
\label{theorem:uniform:upper:bound}
Given $X \in \cT_g$ there exists $C = C(X)$ such that for all $Y \in \cT_g$, 
\begin{displaymath}
F_R(X,Y) \le C(X) e^{h R}.
\end{displaymath}
\end{theorem}

\bold{Notation.} 
\begin{itemize}
\item Let $\mathcal{C}_{g}({\Bbb N})$ be the set of isotopy classes of integral {\it multicurves} on a surface of genus $g.$
\item Given $X, Y \in \te_{g}$ define
$$ M_R(X,Y)= \{\gamma\cdot Y, \gamma \in \Gamma\; |\; \gamma\cdot Y \in B_{R}(X) \},$$ 
so that $ F_{R}(X,Y)= \#(M_R(X,Y)).$
\end{itemize}
\bold{Proof of Theorem~\ref{theorem:volume:growth} assuming
  Theorem~\ref{theorem:uniform:upper:bound}.}
By definition of $F_R$, 
\begin{displaymath}
\mm(B_R(X)) = \int_{\cM_g} F_R(X,Y) \, d\mm(Y) 
\end{displaymath}
We multiply both sides by $e^{-hR}$ and take the limit at $R \to
\infty$. By Theorem~\ref{theorem:uniform:upper:bound} 
we can apply the
bounded convergence theorem to take the limit inside the
integral. Now the theorem follows from
Theorem~\ref{theorem:asympt:count}. 
\qed\medskip

A similar argument yields the following:
\begin{theorem}
\label{theorem:volume:sectors}
Suppose $X \in \cT_g$, and $\cU \subset S(X)$. Then, as $R \to \infty$, 
\begin{displaymath}
\mm(B_R(X) \cap Sect_\cU(X)) \sim \frac{1}{h\, \mm(\cM_g)}e^{h R} \Lambda(Y) \int_\cU
\lambda^-(q) \, ds_X(q).
\end{displaymath}
\end{theorem}

\medskip

\noindent
In the rest of this section we prove
Theorem~\ref{theorem:uniform:upper:bound}. Along the way, we prove 
Theorem~\ref{theorem:Lambda:bounded}. 

\bold{Estimating extremal lengths.}
Consider the Dehn-Thurston parameterization \cite{Harer:Penner:book} of the set of multicurves 
$$ DT: \mathcal{C}_{g}({\Bbb N}) \rightarrow ({\mathbb Z}_{+} \times {\mathbb Z})^{3g-3}$$
defined by 
$$ DT(\beta)=(i(\beta,\alpha_{i}), \tw(\beta,\alpha_{i}))_{i=1}^{3g-3},$$
where $i(\cdot, \cdot)$ denotes the geometric intersection number and
$\tw(\beta,\alpha_{i})$ is the twisting parameter of $\beta$ around $\alpha_{i}$. See $\cite{Harer:Penner:book}$ for more details.
 By a theorem of Bers, we can choose a 
constant $C_g$ depending only on $g$ such that for any surface $Y \in
\cT_g$ there exists a pants decomposition
$\mathcal{P}=\{\alpha_{1},\ldots,\alpha_{3g-3}\}$ 
on $Y$ such that for $1\leq i \leq 3g-3$
$$\Ext_{\alpha_{i}}(Y) \leq C_{g}^{2}.$$ 
We call such a pants decomposition a {\it bounded} pants decomposition for $Y.$
The following result is proved in \cite{Minsky:ext} (see Theorem~5.1,
and equation (4.3)):
\begin{theorem}{\bf (Minsky)}
\label{theorem:minsky:ext}
Suppose $Y \in \cT_g$, and let $\cP = \cP(Y) =
\{\alpha_{1},\ldots,\alpha_{3g-3}\}$ be any bounded pants decomposition on $Y$. 
Then given a simple closed curve $\beta$, $\Ext_\beta(Y)$ 
is bounded from above and below  by  
\begin{equation}\label{es}
\max_{1 \leq j \leq 3g-3}
\left[\frac{i(\beta,\alpha_{j})^{2}}{\Ext_{\alpha_{j}}(Y)}+
\operatorname{tw}^{2}(\beta,\alpha_{j})\Ext_{\alpha_{j}}(Y) \right], 
\end{equation}
up to  a multiplicative constant depending only on $g$.
\end{theorem}

Theorem $\ref{theorem:minsky:ext}$ gives a bound on the  the Dehn-Thurston coordinates of a simple closed
curve in terms of its  extremal length.

\bold{Remark.} The definition of the twist used in equation $(4.3)$ in
\cite{Minsky:ext} is different from the definition we are using
here. We follow the definition used in $\cite{Harer:Penner:book}.$ 
Given a connected simple closed curve $\alpha$ on $\Sigma_g$, let
$h_{\alpha} \in \Gamma$ denote the right Dehn twist around $\alpha$. 
Then in terms of our notation, 
\begin{equation}\label{ours}
\tw(h_{\alpha}^{r}(\beta),\alpha)=\tw(\beta, \alpha)+ r \cdot i(\beta,\alpha).
\end{equation}
 
\begin{coro}\label{esst}
Let $\alpha=\{\alpha_{1},\ldots,\alpha_{3g-3}\}$ be bounded pants decomposition of $X \in \cT_{g}$.
Then there is a constant $c_{1}>0$ such that for any simple closed curve $\beta$ on 
$\Sigma_g$, 
 $$\tw(\beta,\alpha_{i})\leq c_{1} \cdot  \frac{\sqrt{\Ext_{\beta}(X)}}{\sqrt{\Ext_{\alpha_{i}}(X)}} ,$$
 and 
 $$ i(\beta,\alpha_{i})\leq  \sqrt{\Ext_{\beta}(X)}\cdot  \sqrt{\Ext_{\alpha_{i}}(X)} \leq c_{1}  \sqrt{\Ext_{\beta}(X)}.$$
 \end{coro}
\noindent
\bold{Estimating the number of multicurves.} 
Define
$$ E(Y,L )= \#|\{\alpha \in \mathcal{C}_{g}({\Bbb N})\;| \sqrt{\Ext_{\alpha}(Y)} \leq L\}|$$

\noindent Fix $\epsilon_0>0$ small enough such that if $\alpha$ and
$\beta$ on $Y \in \cT_{g}$ satisfy $\Ext_{\alpha}(Y) \leq
\epsilon_0^{2}$ and $\Ext_{\beta}(Y) \leq C_{g}^{2}$ then $i(\alpha,
\beta)=0.$ Note that any bounded pants decomposition $\mathcal{P}$ on
$Y$ includes all simple closed curves of extremal length $ \leq
\epsilon_0^2$ on $Y.$
 
 Let
$$G(Y) = 1+ \prod_{\Ext_\gamma(Y) \le \epsilon_0^2}
\frac{1}{\sqrt{\Ext_\gamma(Y)}},$$ 
and the product ranges over all
simple closed curves $\gamma$ on the surface $Y$ with $\Ext_\gamma(Y)
\le \epsilon_0^2$.

Using Theorem~\ref{theorem:minsky:ext}, we obtain the following:

\begin{theorem}\label{count}
There exists a constant $C>0$ such that for every $Y \in \te_{g}$ and $L>0$ we have
$$ E(Y,L) \leq C \cdot G(Y) \cdot L^{6g-6}.$$ 
Moreover, for any $Y$, if $1/L$ is bounded by an absolute constant
times the square root of the extremal length of the shortest curve on
$Y$, then
$$E(Y,L) \leq C L^{6g-6}.$$
\end{theorem}
\noindent
{\bf Sketch of the proof.}
In order to use the bound given by equation ($\ref{es}$), first we 
fix a bounded pants decomposition $$\mathcal{P}=\{\alpha_{1},\ldots,\alpha_{3g-3}\}$$ on $Y$. 
Note that, this pants decomposition should include all small closed curves on $Y.$ 
Let $m_{i}=i(\beta,\alpha_{i})$, $t_{i}=\tw(\beta,\alpha_{i})$, and $s_{i}=\sqrt{\Ext_{\alpha_{i}}(Y)}$. 
So by equation ($\ref{es}$), $\Ext_{\beta}(Y) \leq L^{2}$ implies that  
$\frac{m_{i}}{s_{i}}+ |t_{i}| s_{i} = O(L)$.
We use the following elementary lemma:
\begin{lemma}
For $s>0$, define $A_{s}(L)$ by
$$A_{s}(L)=\{(a,b)\;|\; a,b \in {\mathbb Z}_{+}, a \cdot s+\frac{b}{s} \leq L\} \subset {\mathbb Z}_{+} \times {\mathbb Z}_{+}.$$
Then for any $L>0$, $|A_{s}(L)| \leq 4 \max\{s,1/s\}\cdot L^{2}.$
For $L> \max\{s,\frac{1}{s}\}$, we have $A_{s}(L) \leq 4 L^{2}.$ 
\end{lemma}

Now applying the preceding lemma to the $s_{i}$, and using the
Dehn-Thurston parameterization of multicurves, we get 
$$E(Y,L) \leq \prod\limits_{i=1}^{3g-3} |A_{s_{i}}(L)| \leq C\; G(Y) L^{6g-6}.$$
\hfill $\Box$

\bold{Proof of Theorem~\ref{theorem:Lambda:bounded}.}
This follows from the second part of Theorem~\ref{count}.
Recall (see \S\ref{sec:notation}) that the extremal length can be extended continuously to a map  $$\Ext: 
\mf \times \te_{g} \rightarrow {\mathbb R}_{+}$$ such that 
$\Ext_{t \lambda}(X) =t^{2} \Ext_{\lambda}(X).$
Also the space $\mf$ has a piecewise linear integral structure, and elements of $\cS_{g}$ are in one to one correspondence with the integral points.
Hence, 
\begin{equation}
\label{eq:Lambda:and:multicurves}
\Lambda(X)=\vo\{\eta \in \mf \; | \sqrt{\Ext_{\eta}(X)} \leq 1\}=
\frac{\vo\{\eta \in \mf \; | \sqrt{\Ext_{\eta}(X)} \leq L\}}{L^{6g-6}} =
\lim_{L \rightarrow \infty} \frac{E(x,L)}{L^{6g-6}},
\end{equation}
where to justify the last equality use the fact that for any as $L \to \infty$ the number of lattice
points in the dilation $LA$ is asymptotic to the volume of $LA$, where $A = \{\eta \in \mf \; | \sqrt{\Ext_{\eta}(X)} \leq 1\}$.

On the other hand, by the second assertion in 
Theorem~\ref{count} when $L$ is big enough
$E(Y,L) \leq C L^{6g-6}$. 
Therefore, $\Lambda(Y) \leq C,$ where the bound does not depend on
$Y.$
\qed\medskip

\bold{Remark.} Let $p: \te_{g} \to \mathcal{M}_{g}$ be the natural
projection. Fix a compact subset $\cK \subset \mathcal{M}_{g}$. For
simplicity, we will also denote $p^{-1}(\cK)\subset \te_{g}$ by $\cK$.
If  
$\mathcal{P}=\{\alpha_{i}\}$ is a bounded pants decomposition on 
$X\in \cK$ then for any $1 \leq i \leq 3g-3$,
$$\frac{1}{C_{\cK}} \leq \sqrt{\Ext_{\alpha_{i}}(X)} \leq C_{\cK},$$ where
$C_{\cK}$ is a constant which only depends on $\cK.$ 
Therefore, by Corollary~\ref{esst}, for the elements in the thick part of the
Teichm\"uller space, the Dehn-Thurston coordinates of a simple closed
curve are bounded by its extremal length.
More precisely, if $\alpha=\{\alpha_{1},\ldots,\alpha_{3g-3}\}$ be
bounded  pants decomposition of $X \in \cK$. 
then there is a constant $c_{1}>0$, depending only on $\cK$, such that
for any simple closed curve $\beta$ on  
$\Sigma_g$, 
 $$\tw(\beta,\alpha_{i})\leq c_{1} \cdot  \sqrt{\Ext_{\beta}(X)}, \qquad i(\beta,\alpha_{i})\leq c_{1}\cdot  \sqrt{\Ext_{\beta}(X)}.$$
 
\noindent
\begin{proposition}
\label{prop:multicurves}
There exists a constant $C_{3}$ depending only on $\cK$ such that for every 
$X \in \cK$, $Y\in \te_{g}$ and $R>0$ we have  
$$F_{R}(X,Y) \leq C_{3}\; E(Y, e^{R}).$$
\end{proposition}

Our goal is to assign to any point $Z= \gamma \cdot Y \in M_R(X,Y)$ a
unique integral multicurve $\beta_{Z} \in E(Y,e^{R}).$ This would
imply that $$F_{R}(X,Y) \leq C e^{(6g-6) R}.$$ (Note that we are
assuming that $X$ is in the compact part of $\cT_g$ and $d(X,Y) < R$;
this implies that the shortest curve on $Y$ has extremal length at
least constant times $e^{-2 R}$, hence we may use the second statement
of Theorem~\ref{count}.)

The most natural candidate for $\beta_{Z}$ is $\gamma^{-1} \alpha$
where $\alpha$ is a fixed pants decomposition of $X$. However this
correspondence is not one to one. Therefore, we need to modify the
construction.

Given $Z=\gamma \cdot Y \in M_R(X,Y)$, and $\alpha \in \cS_{g}$ we have 
$$\frac{\sqrt{\Ext_{\alpha}(Z)}}{\sqrt{\Ext_{\alpha}(X)}} \leq e^{R}.$$
So by setting $L=e^{R}$, and $\beta=\gamma^{-1}(\alpha)$ we have
\begin{equation}\label{onee}
\sqrt{\Ext_{\beta}(Y)} \leq L \cdot \sqrt{\Ext_{\alpha}(X)}.
\end{equation}

Given $Z=\gamma \cdot Y \in M_R(X,Y)$, let  $\alpha(Z)_{i}=\gamma^{-1} \alpha_{i} \in \cS_{g}.$
Then from ($\ref{onee}$) we get
$$\sqrt{\Ext_{\alpha(Z)_{i}}(Y)} \leq C_{\cK} \cdot L .$$

Moreover, for  $Z_{1}=\gamma_{1} Y,Z_{2}=\gamma_{2} Y \in \Gamma\cdot
Y$, we have $\alpha(Z_{1})=\alpha(Z_{2})$ if and only if there are
$r_{1},\ldots r_{3g-3} \in{\mathbb Z}$ such that
$$\gamma_{1}=h_{\alpha_{1}}^{r_{1}} \cdots h_{\alpha_{3g-3}}^{r_{3g-3}}
\cdot \gamma_{2}.$$  
Moreover, we have:
\begin{lemma}\label{Mmain}
Let $\alpha=\{\alpha_{1},\ldots,\alpha_{3g-3}\}$ be a bounded pants
decomposition on $X$ and suppose $Y_{0} \in B_{R}(X)$.
If $h_{\alpha_{1}}^{r_{1}} \cdots h_{\alpha_{3g-3}}^{r_{3g-3}} (Y_{0}) \in B_{R}(X)$ then for $1 \leq i\leq 3g-3$
$$|r_{i}| \cdot \sqrt{\Ext_{\alpha_{i}}(Y_{0})} \leq C_{2} \frac{e^R} {\sqrt{\Ext_{\alpha_{i}}(X)}}.$$
Here $C_{2}$ is a constant independent of $X$ and $Y_{0}.$

\end{lemma}
\noindent
{\bf Sketch of the proof of Lemma $\ref{Mmain}$.}
Let $s_{i}=\sqrt{\Ext_{\alpha_{i}}(Y_{0})},$ and $L=e^{R}$.
The key point in the proof is Lemma 6.3 in \cite{Minsky:ext}. From this lemma, for any $i$,
there exists a multicurve $\beta_{i}$ of extremal length $b_{i}^2=\Ext_{\beta_{i}}(Y_{0})$ such that  
$$i(\beta_{i},\alpha_{i}) \geq c s_{i}\cdot b_{i}.$$
Since $d(X,Y_{0}) \leq R$, $\sqrt{\Ext_{\beta_{i}}(X)} \leq b_{i}\cdot L.$ Also the assumption 
$$d(h_{\alpha_{1}}^{-r_{1}}\cdots h_{\alpha_{3g-3}}^{-r_{3g-3}}(X), Y_{0}) \leq R$$ 
implies that 
\begin{equation}\label{sec}
\sqrt{\Ext_{h_{\alpha_{1}}^{r_{1}}\cdots h_{\alpha_{3g-3}}^{r_{3g-3}}(\beta_{i})}(X))} =\sqrt{\Ext_{\beta_{i}}(h_{\alpha_{1}}^{-r_{1}}\cdots h_{\alpha_{3g-3}}^{-r_{3g-3}}(X))} \leq b_{i} \cdot L.
\end{equation}
Then 
\begin{enumerate}

\item Since  $\sqrt{\Ext_{\beta_{i}}(X)} \leq b_{i} \cdot L,$ Corollary \ref{esst} implies that the twist and intersection coordinates of the curve $\beta_{i}$ are bounded by a multiple of $b_{i}\cdot L$; in particular $$|\tw(\beta_{i},\alpha_{i})|=O\left( \frac{ b_{i}\cdot L} {\sqrt{\Ext_{\alpha_{i}}(X)}}\right).$$
\item Similarly, by equation ($\ref{sec}$), applying Corollary \ref{esst} for $h_{\alpha_{1}}^{r_{1}}\cdots h_{\alpha_{3g-3}}^{r_{3g-3}}(\beta_{i})$ on $X$ 
with respect to $\{\alpha_{i}\}$ implies that 
$$ |\tw(h_{\alpha_{1}}^{r_{1}}\cdots h_{\alpha_{3g-3}}^{r_{3g-3}}(\beta_{i}), \alpha_{i})|= O\left( \frac{ b_{i}\cdot L} {\sqrt{\Ext_{\alpha_{i}}(X)}}\right).$$
\item On the other hand, by the definition (see equation ($\ref{ours}$)) we have $$|\tw(h_{\alpha_{1}}^{r_{1}}\cdots h_{\alpha_{3g-3}}^{r_{3g-3}}(\beta_{i}), \alpha_{i})|=|r_{i} \cdot i(\alpha_{i},\beta_{i})+\tw(\beta_{i},\alpha_{i})| \geq $$ 
$$ \geq |r_{i}|\cdot i(\alpha_{j},\beta_{i})-|\tw(\beta_{i},\alpha_{i})|.$$
\end{enumerate}
 So we have 
$$|r_{i}|\cdot s_{i}\cdot b_{i} \leq \frac{1}{c}\; |r_{i}| \cdot i(\beta_{i},\alpha_{i}) \leq C_{2} L\cdot b_{i} \Longrightarrow  |r_{i}|\cdot s_{i} \leq \frac{C_{2} L} {\sqrt{\Ext_{\alpha_{i}}(X)}}.$$
\hfill $\Box$\\
\noindent
{\bf Remark.}
As a result, if the assumption of Lemma $\ref{Mmain}$ holds and $X \in \cK$ then 
$$|r_{i}| \cdot \sqrt{\Ext_{\alpha_{i}}(Y_{0})} \leq C e^R,$$
where $C$ is a constant which only depends on $\cK.$

We remark that for any two disjoint simple closed curves $\eta_{1}$
and  $\eta_{2}$, and $m_{1},m_{2}>0$, we have, by the definition of
extremal length, 
$$\sqrt{\Ext_{m_{1}\cdot \eta_{1}}(X)}+ \sqrt{\Ext_{m_{2} \cdot \eta_{2}}(X)} \geq \sqrt{\Ext_{m_{1}\cdot \eta_{1}+m_{2} \cdot \eta_{2}}(X)}.$$
\begin{coro}\label{corol}
If the assumption of Lemma $\ref{Mmain}$ holds and $X \in K$ then for $\widehat{\alpha}= \sum\limits_{i=1}^{3g-3} |r_{i}|\cdot  \alpha_{i}$ we have 
$$\sqrt{\Ext_{\widehat{\alpha}}(Y_{0})} \leq C_{3} e^R.$$
\end{coro}

Therefore, $\widehat{\alpha}$ also defines a multicurve of extremal
length bounded by $e^{2 R}$. This completes the proof of 
Proposition~\ref{prop:multicurves}, and
thus in view of Theorem~\ref{count} the proof of
Theorem~\ref{theorem:uniform:upper:bound}. 
\qed\medskip

\appendixmode

\section{Appendix: Proof of Theorem~\ref{theorem:sx:multiple:zero:locus}.}

In this section, $\qstartg$ denotes the space of non-zero quadratic
differentials on marked compact surfaces of genus $g$, without
restriction on area. 

\bold{A local coordinate system near the multiple zero locus.}
Near the multiple zero locus, the coordinate system given by the
period map (see Lemma~\ref{lemma:good:basis}) is singular. Instead we use an
alternative coordinate system from \cite{hubmas}.

Suppose $q_0 \in \qstartg$ has multiple zeroes, say $w_1, \dots, w_m$. 
Let $m_i$ be the multiplicity of $q_0$ at $w_i$. Then, by 
\cite[Proposition~3.1]{hubmas}, for $1 \le i \le m$ 
there exists a choice of local coordinate $z_i$, mapping $w_i$ to $0$,
such that any $q \in \qstartg$ near $q_0$, on a neighborhood of $w_i$, has
the form:
\begin{displaymath}
q = \left(z_i^{m_i} + \sum_{j=0}^{m_i-2} a_{ij} z_i^j\right)
(dz_i)^2. 
\end{displaymath}

The coordinates $z_i$ are uniquely determined by $q$. 
Fix some $\delta > 0$, and let
$w_i'$ be the point corresponding to $z_i = \delta$. For $q \in \qstartg$
near $q_0$ let $\tilde{X}$ denote the canonical double cover which makes the
foliations corresponding to $q$ orientable, and let $\sigma$ be the
involution of $\tilde{X}$ so that $\tilde{X}/\sigma$ is the
surface $X$ (with the flat structure given by $q$). Let $\Sigma_\delta
= \sigma^{-1}(w_1', \dots, w_m')$, and let $H$ denote the relative homology
group $H_1(\tilde{X}, \Sigma_\delta, \zed)$. Let $H_{odd}$ denote the
odd part of $H$ under the action of the involution $\sigma$. Let
$\gamma_1, \dots, \gamma_m$ be a integral basis for $H_{odd}$. 
Let 
\begin{displaymath}
\lambda_i = \int_{\gamma_i} \sqrt{q}. 
\end{displaymath}


Suppose $m_i$ is even. Let $w_{ij}$, $1 \le j \le m_i$ denote the
zeroes of $q$ which tend to $w_i$ as $q \to q_0$. 
Let $\eta_i$ be a fixed small circle in the coordinates $z_i$, 
centered at $z_i = 0$ (i.e. $w_i$). 
Then  if $q$ is sufficiently close to $q_0$, 
$\eta_i$ separates all the $w_{ij}$ from the rest of the surface. 
Let $b_i = \int_{\eta_i} \sqrt{q}$. Then, a residue calculation shows
that 
\begin{equation}
\label{eq:middle:coeff}
b_i = a_{\frac{m_i}{2}-1} + \text{some polynomial in $a_{ij}$, $m_i -2
  \ge j \ge   \frac{m_i}{2}$. }
\end{equation}

\begin{proposition}[Hubbard-Masur]
\label{prop:hubbard:masur:coordinates}
If $q_0$ is not a square of an Abelian differential, then there exists
a neighborhood $U \subset \qstartg$ of $q_0$  such that 
the $a_{ij}$ and the $\lambda_i$ are smooth local coordinates on $U$. 
If $q_0$ is the square of an Abelian differential, then there exists a
neighborhood $U$ of $q_0$ such that under the constraint $\sum_i b_i =
0$, the $a_{ij}$ and the $\lambda_i$ are local coordinates on $U$. 

In both cases, near $q_0$, the natural projection map $\pi:
\qstartg \to \cT_g$ is a submersion in these coordinates. 
\end{proposition}

\bold{Proof.} See \cite[Proposition 4.7]{hubmas}. 
\qed\medskip

The constraint
$\sum_i b_i = 0$ appears because in the case when $q_0$ is the square
of an Abelian differential, the part of the surface which is outside 
all the circles $\eta_i$ has an orientable foliation. In that case, in
view of (\ref{eq:middle:coeff}), we
may drop $a_{\frac{m_i}{2}-1}$ for some $i$ from the coordinate system.

We introduce a covering of $\qstartg$ by open sets
$U_\alpha$, such that  each $U_\alpha$ is the neighborhood $U$ of 
Proposition~\ref{prop:hubbard:masur:coordinates}, for some $q_0 \in
\qstartg$. We may assume that the $U_\alpha$ are invariant under the
operation of multiplying the quadratic differential by a real
number. We also assume that covering is uniformly locally finite
on compact sets, i.e. that for any compact set $\cK \in \cT_g$, there
exist a number $N$ depending only on $\cK$ such that each point in
$\hat{\pi}^{-1}(\cK)$ belongs to at most $N$ sets $U_\alpha$. 

Let $\psi_\alpha$ be a partition of unity subordinate to $U_\alpha$. 
For $q, q' \in U_\alpha$, let $\{a_{ij}$, $\lambda_i\}$ be the coordinates
of $q$, and $\{a_{ij}'$, $\lambda_i'\}$ be the coordinates of $q'$. 
Let $D_\alpha(q,q') = \sum_{ij} |a_{ij} - a_{ij}'| + \sum_{i}
|\lambda_i - \lambda_i'|$. For $q$ and $q'$ sufficiently close so 
that they belong to the same $U_\alpha$, let $D(q,q') = \sum_\alpha \psi_\alpha
D_\alpha(q,q')$. 

Recall $d_{\cT}(X,Y)$ denotes the Teichm\"uller distance between $X$
and $Y$, and $S(X)$ denotes the sphere of unit area quadratic
differentials which are holomorphic at $X$. 
\begin{lemma}
\label{lemma:coeffs:vs:teich:distance}
Suppose $X$, $Y$ in $\cT_g$ are
sufficiently close and belong to a compact set $\cK$. 
Then, there exist $c > 1$ depending only on $\cK$ such that 
\begin{equation}
\label{eq:coeffs:to:Teich}
c^{-1} \inf_{q \in S(X)} D(q, \pi^{-1}(Y))  \le d_{\cT}(X,Y) \le
c \sup_{q \in S(X)} D(q, \pi^{-1}(Y)). 
\end{equation}
\end{lemma}

\bold{Proof.} Let $K$ be the pullback of $\cK$ from $\cT_g$ to
$\qonetg$. Then $K$ is also compact. Near any $q \in K$, by 
Proposition~\ref{prop:hubbard:masur:coordinates}, the map $\pi: \qstartg
\to \cT_g$ is a submersion, and therefore can be approximated by 
a linear map, for which (\ref{eq:coeffs:to:Teich}) is immediate. 
The rest follows by compactness. 
\qed\medskip

Let $\beta_\alpha$ be the measure on $U_\alpha \subset\qstartg$ 
given by $\prod_{i,j} da_{ij} \prod_{i} d\lambda_i$. Let $\beta =
\sum_\alpha \psi_\alpha \beta_\alpha$.

For $V \subset \pi^{-1}(X)$, let $V^* \subset \qstartg$ denote the 
set of all quadratic differentials which have the same horizontal
foliation as some $q \in V$. 
\begin{lemma}
\label{lemma:small:measure:in:coeff:coordinates}
For every $\delta > 0$, any $X \in \cT_g$  and any $q_0 \in \qstartg$ there
exists an  open subset $V$ of $\pi^{-1}(X)$ containing $q_0$, such that 
for all sufficiently small $\epsilon > 0$, 
\begin{displaymath}
\beta(V^* \cap \pi^{-1}(B(X,\epsilon))) < \delta \epsilon^{6g-6}. 
\end{displaymath}
\end{lemma}

\bold{Proof.} Fix $\delta' > 0$. Recall that by a horosphere (e.g. leaf
of $\cF^{ss}$) in
$\qstartg$ we mean the set of quadratic differentials with a fixed 
horizontal foliation. Since the foliation by horospheres is continuous, 
there exists a neighborhood $V$ of $q_0 \in \pi^{-1}(X)$ such that
$V \subset B_D(q_0, \delta')$, where $B_D$ denotes a ball in the
metric $D$. Then, for $\epsilon > 0$ sufficiently small, and in view of
Lemma~\ref{lemma:coeffs:vs:teich:distance}, there exist constants
$C_1$, $C_2$ such that 
\begin{equation}
\label{eq:contained:in:hbhd:ball}
\pi^{-1}(V^* \cap B(X,\epsilon))  \subset \{ q \in \qstartg \st D\left(q,
B_D(q_0, C_2 \delta') \cap \pi^{-1}(X)\right)  < C_1 \epsilon \} . 
\end{equation}
Since $\pi^{-1}(X)$ is smooth, the $\beta$-measure of the set on the
right of (\ref{eq:contained:in:hbhd:ball}) is $\epsilon^{6g-6}
c(\delta')$, where $c(\delta') \to 0$ as $\delta' \to 0$. This implies
the lemma.
\qed\medskip

\begin{proposition}
\label{prop:jacobian}
The measure $\mu$ on $\qstartg$ is absolutely continuous with respect
to $\beta$, and for any compact subset $\cK$ of $\cT_g$ there exists a
constant $C$ such that for all $q \in \pi^{-1}(\cK)$, 
$\left|\frac{d\mu}{d\beta}(q)\right| < C$. Also, away from the multiple
zero locus, $\left|\frac{d\mu}{d\beta}(q)\right|$ is a smooth
non-vanishing function of $q$.
\end{proposition}

\bold{Remark.} It is also possible to show that there exists a
constant $C'$ depending only on $\cK$ such that for all $q \in
\pi^{-1}(\cK)$, one has $C' < \left|\frac{d\mu}{d\beta}(q)\right|$. Since we do
not need this, we will omit the proof. 

\bold{Proof of Theorem~\ref{theorem:sx:multiple:zero:locus}, assuming
Proposition~\ref{prop:jacobian}.}
By Lemma~\ref{lemma:coeffs:vs:teich:distance},
$\beta(\pi^{-1}(B(X,\epsilon))) = O(\epsilon^{6g-6})$. Therefore, by
Proposition~\ref{prop:jacobian}, there exist $0 < c_1 < c_2$ such that
\begin{equation}
\label{eq:equation:coeffs:ball}
c_1 \epsilon^{6g-6} < \mu\left(\pi^{-1}(B(X,\epsilon))\right) < c_2
\epsilon^{6g-6}. 
\end{equation}
Now, 
\begin{displaymath}
ds_X(V) = \lim_{\epsilon \to 0} \frac{\mu(V^* \cap
  \pi^{-1}(B(X,\epsilon)))}{\mu(\pi^{-1}(B(X,\epsilon)))} \le
\limsup_{\epsilon \to 0} \frac{C \beta(V^* \cap
  \pi^{-1}(B(X,\epsilon)))}{\mu(\pi^{-1}(B(X,\epsilon)))} \le \frac{C
  \delta \epsilon^{6g-6}}{c_1 \epsilon^{6g-6}}.
\end{displaymath}
where for the first inequality we used
Proposition~\ref{prop:jacobian}, and for the last estimate we used 
Lemma~\ref{lemma:small:measure:in:coeff:coordinates} and
(\ref{eq:equation:coeffs:ball}). Since $\delta$ is arbitrary, the
theorem follows. 
\qed\medskip

The rest of the section will consist of the proof of
Proposition~\ref{prop:jacobian}. To simplify notation, we work
with one zero of $q_0$ at a time. 
We write:
\begin{displaymath}
z^{m} + \sum_{i=0}^{m-1} a_i z^i = \prod_{j=1}^m (z - z_j). 
\end{displaymath}
Then the $a_i$ are symmetric polynomials in the $z_j$. It is well
known that the Jacobian of the map from the zeroes $z_j$'s to the
coefficients $a_i$'s
is the Vandermonde determinant. 
In other words, if we use the notation $\frac{\partial(y_1,\dots,
  y_n)}{\partial(x_1,\dots, x_n)}$ for the Jacobian determinant of the matrix
$\left\{\frac{\partial y_i}{\partial x_j}\right\}_{1\le i,j \le n}$, we have:
\begin{equation}
\label{eq:jac:symmetric:functions}
\frac{\partial(a_0, \dots, a_{m-1})}{\partial(z_1, \dots, z_m)} =
\prod_{i < j} (z_i - z_j). 
\end{equation}
(This follows from comparing degrees and noting the anti-symmetry of
the determinant).

We choose a basis for the homology relative to the zeroes. This
amounts to choosing a spanning subtree $T$ from the complete graph
connecting the zeroes $z_j$. We may choose the tree
in such a way that for any $z_i$ and $z_j$, $|z_i - z_j|$ is within a 
multiplicative constant 
(depending only on $m$) of the length of
the path in $T$ connecting $z_i$ and $z_j$. 
Let the (oriented) edges of $T$ be $e_1,
\dots, e_{m-1}$. Let $e_k^+$ be the zero at the head of $e_k$, and
$e_k^-$ the zero at the tail of $e_k$. We write $\vec{e}_k = e_k^{+} -
e_k^{-}$ (so $\vec{e}_k \in \cx$).

In our setting, we need to restrict to $a_{m-1} = 0$. 
Note that $a_{m-1} = z_1 + \dots + z_m$. It follows that
\begin{equation}
\label{eq:jac:part:coeffs:to:zeroes}
\frac{\partial(a_0, \dots, a_{m-2})}{\partial(\vec{e}_1, \dots,
  \vec{e}_{m-1})} = \frac{\partial(a_0, \dots,
  a_{m-1})}{\partial(\vec{e}_1, \dots,  \vec{e}_{m-1}, a_{m-1})} = 
\frac{\partial(a_0, \dots, a_{m-1})}{\partial(z_1, \dots,
  z_m)} \frac{\partial(z_1, \dots, z_m)}{\partial(\vec{e}_1, \dots,
  \vec{e}_{m-1}, a_{m-1})} =  \prod_{i < j} (z_i - z_j),
\end{equation}
where we have used (\ref{eq:jac:symmetric:functions}). 
 
If $z_j$ is a zero, and $e$ is an edge of $T$, let 
\begin{displaymath}
d_+(z_j,e) = \max(|z_j-e^-|, |z_j - e^+|). 
\end{displaymath}
We will need the following combinatorial lemma:
\begin{lemma}
\label{lemma:strange:comb}
\begin{equation}
\label{eq:strange:estimate}
\prod_{e \in T} \prod_{p=1}^m d_+(z_p, e)^{1/2} \le C \prod_{i < j} |z_i
- z_j|,
\end{equation}
where $C$ depends only on $m$. 
\end{lemma}
\bold{Proof of Lemma~\ref{lemma:strange:comb}.}
Let $e'$ be the longest edge of
$T$. If we cut along $e'$, we separate the tree into two subtrees
$T_1$ and $T_2$
say of size $m_1$ and $m_2$. If $z_i \in T_1$ and $z_j \in T_2$, then
by the assumption on $T$, $|z_i - z_j|$ is comparable to the length
of the path in the tree connecting $z_i$ to $z_j$. This path contains
$e'$, which is by assumption the longest edge in the tree. Therefore, 
if $z_i \in T_1$ and $z_j \in T_2$ then $|z_i - z_j|$ is within a 
multiplicative constant of $e'$. 
Hence, the left hand side of (\ref{eq:strange:estimate})
is within a multiplicative constant of 
\begin{displaymath}
|e'|^{m_1 m_2} \prod_{\stackrel{i < j}{z_i, z_j \in T_1}} |z_i - z_j | 
\prod_{\stackrel{i < j}{z_i, z_j \in T_2}} |z_i - z_j |. 
\end{displaymath}
To estimate the right hand side of (\ref{eq:strange:estimate}), 
note that each factor of the form $d_+(z_p, e)^{1/2}$ where $z_p \in T_1$
and either $e \in T_2$ or $e = e'$ is within a multiplicative constant of 
$|e'|^{1/2}$. The number of such factors is $m_1 m_2$. Also the same
factors appear when $z_p \in T_2$ and $e \in T_1$ or $e = e'$. Then
the right hand side of (\ref{eq:strange:estimate}) is within a
multiplicative constant of 
\begin{displaymath}
|e'|^{m_1 m_2} \prod_{e \in T_1} \prod_{z_p \in T_1} d_+(z_p, e)^{1/2}
\prod_{e \in T_2} \prod_{z_p \in T_2} d_+(z_p, e)^{1/2}
\end{displaymath}
The estimate (\ref{eq:strange:estimate}) now follows by induction. 
\qed\medskip

\begin{lemma}
\label{lemma:periods:to:zeros}
Suppose $z_1, \dots, z_m$ are in some bounded set $\cK$. 
Let 
\begin{displaymath}
\Omega_k = \int_{e_k^-}^{e_k^+} \sqrt{(z-z_1) \dots (z-z_n)} \, dz,
\end{displaymath}
so that $\Omega_k$ is the holonomy of the edge $e_k$. 
Then,
\begin{displaymath}
\left| \frac{\partial(\Omega_1, \dots,
    \Omega_{n-1})}{\partial(\vec{e}_1, \dots, \vec{e}_{m-1})} \right|
\le C \prod_{i < j} |z_i - z_j|,
\end{displaymath}
where $C$ depends only on $m$ and $\cK$. 
\end{lemma}

\bold{Proof.} Suppose for the moment that $z_1, \dots, z_m$ are
independent variables. 
We claim that 
\begin{equation}
\label{eq:derivative:of:period}
\left|\frac{\partial \Omega_k}{\partial z_j} \right| \le C
\prod_{p=1}^n d_+(z_p,e_k)^{1/2},
\end{equation}
where $C$ depends only on $m$ and $\cK$. Indeed, 
\begin{displaymath}
\frac{\partial \Omega_k}{\partial z_j} =  
\int_{e_k^-}^{e_k^+} \frac{ \sqrt{(z-z_1) \dots (z-z_n)}}{2
  (z-z_j)}  \, dz
\end{displaymath}
If we write the 
numerator in the integral as $P (z-z_j)^{1/2}$, then $P$ 
is bounded by a constant times
$\prod\limits_{\stackrel{p=1}{p\ne j}}^n d_+(z_p,e_k)^{1/2}$. Then, 
the integral is bounded by
\begin{displaymath}
\frac{|P|}{2} \left|\int_{e_k^-}^{e_k^+} \frac{dz}{(z-z_j)^{1/2}}\right|.
\end{displaymath}
If $d_+(z_j,e_k) > \frac{1}{2} |e_k^+-e_k^-|$, 
\begin{displaymath}
\left|\int_{e_k^-}^{e_k^+} \frac{dz}{(z-z_j)^{1/2}}\right| \le C
\frac{|e_k^+ - e_k^-|}{d_+(z_j,e_k)^{1/2}} \le C' |e_k^+-e_k^-|^{1/2}
\end{displaymath}
(where we have used the assumption that $|e_k^+-e_k^-|$ is bounded). 
If $d_+(z_j,e_k) < \frac{1}{2} |e_k^+-e_k^-|$, 
then 
\begin{displaymath}
\left|\int_{e_k^-}^{e_k^+} \frac{dz}{(z-z_j)^{1/2}}\right| \le C''
|e_k^+-e_k^-|^{1/2},
\end{displaymath}
where $C''$ depends only on $m$ and $\cK$. 
Thus, 
\begin{displaymath}
\left| \frac{\partial \Omega_k}{\partial z_j} \right| \le C
|e_k^+-e_k^-|^{1/2} \prod_{p \ne j} d_+(z_p, e_k)^{1/2}. 
\end{displaymath}
But by the choice of the tree $T$, \mc{explain} for any $j$, 
$d_+(z_j, e_k) \ge c |e_k^+-e_k^-|$, where $c$ depends only on $m$. 
Thus, (\ref{eq:derivative:of:period}) holds. 

Since $\sum_{j=1}^m z_j = 0$, there are $m-1$ linearly independent
$z_j$'s and therefore we can express $z_j =
\sum_{k=1}^{m-1} \eta_{jk} \vec{e}_k$, where the $\eta_{jk}$ are
bounded depending only on $m$. Then, we get
\begin{displaymath}
\left|\frac{\partial \Omega_j}{\partial \vec{e}_k}\right| \le C'
\prod_{p=1}^m d_+(z_p, e_k)^{1/2}. 
\end{displaymath}
Lemma~\ref{lemma:periods:to:zeros} now follows from
Lemma~\ref{lemma:strange:comb}. 
\qed\medskip

\bold{Proof of Proposition~\ref{prop:jacobian}.} It is enough to show
that for any $\alpha$, $|\frac{\partial \mu}{\partial \beta_\alpha}|
\le C$ on $U_\alpha$. As above, we work with one zero of $q_0$ at a
time. By (\ref{eq:jac:part:coeffs:to:zeroes}) and
Lemma~\ref{lemma:periods:to:zeros}, 
\begin{equation}
\label{eq:jac:almost:done}
\left| \frac{\partial(\Omega_1, \dots,
    \Omega_{m-1})}{\partial(a_0, \dots, a_{m-2})} \right|
= 
\left| \frac{\partial(\Omega_1, \dots,
    \Omega_{m-1})}{\partial(\vec{e}_1, \dots, \vec{e}_{m-1})} \right|
\left| \frac{\partial(\vec{e}_1, \dots,
    \vec{e}_{m-1})}{\partial(a_0, \dots, a_{m-2})} \right| \le C,
\end{equation}
where $C$ depends only on $m$ and $\cK$. 

Let $\mu' = d\Omega_1 \, d\bar{\Omega}_1 \dots  d\Omega_{m-1} \,
d\bar{\Omega}_{m-1}$, and let $\beta' = da_0 \, 
d\bar{a}_0  \dots da_{m-2} \, d\bar{a}_{m-2}$. Then, by
(\ref{eq:jac:almost:done}), 
\begin{displaymath}
\frac{d\mu'}{d\beta'} = \left| \frac{\partial(\Omega_1, \dots,
    \Omega_{m-1})}{\partial(a_0, \dots, a_{m-2})} \right|^2 \le C^2. 
\end{displaymath}
Now, 
\begin{displaymath}
\frac{d\mu}{d\beta_\alpha} = \prod_j \frac{d\mu'_j}{d\beta'_j},
\end{displaymath}
where the product is over distinct zeroes of $q_0$. Now
Proposition~\ref{prop:jacobian} follows from
(\ref{eq:jac:almost:done}). 
\qed

\bigskip
Department of Mathematics, University of Illinois, Urbana, IL 61801 USA;\\
jathreya@illinois.edu

\bigskip
Department of Mathematics,
Rice University, MS 136, 6100 Main Street, Houston, Texas
77251-1892 \textrm{and} The Steklov Institute of
Mathematics, Russian Academy of Sciences, Gubkina str. 8,
119991, Moscow, Russia; \\
 aib1@rice.edu \textrm{and} bufetov@mi.ras.ru

\bigskip
Department of Mathematics, University of Chicago, Chicago IL~60637,
USA; \\
eskin@math.uchicago.edu

\bigskip
Department of Mathematics, Stanford University, Stanford CA 94305 USA;\\
mmirzakh@math.stanford.edu

\end{document}